\documentclass[AMS,STIX1COL]{WileyNJD-v2}
\articletype{Survey}%
\usepackage{cleveref}

\def\vek#1{\mathbf{#1}}
\def\bbC{\mathbb{C}} 

\def\R{\mathbb{R}}
\def\K{\mathbb{K}}
\def\CS{{\cal S}}

\def\CC{{\cal C}}
\def\CK{{\cal K}}

\def\CU{{\cal U}}
\def\CV{{\cal V}}

\def\CX{{\cal X}}
\def\CY{{\cal Y}}

\newcommand{\argmin}[1]{\underset{#1}{\text{{\rm argmin}}}}
\newcommand{\nn}{\nonumber}

\def\mand{\mbox{\ \ \ and\ \ \ }}

\def\mwith{\mbox{\ \ \ with\ \ \ }}
\def\mwhere{\mbox{\ \ \ where\ \ \ }}

\def\msuchthat{\mbox{\ \ \ such that\ \ }}

\def\bq{\begin{quotation}}
\def\eq{\end{quotation}}

\def\b{\beta}

\def\z{\zeta}

\def\P{\Pi}

\def\t{\tau}

\newcommand{\M}[1]{ \mathbf{#1} }  

\newcommand{\MA}{\M{A}}
\newcommand{\MB}{\M{B}}
\newcommand{\MC}{\M{C}}
\newcommand{\MD}{\M{D}}
\newcommand{\ME}{\M{E}}

\newcommand{\MG}{\M{G}}
\newcommand{\MH}{\M{H}}
\newcommand{\MI}{\M{I}}

\newcommand{\ML}{\M{L}}

\newcommand{\MN}{\M{N}}

\newcommand{\MP}{\M{P}}
\newcommand{\MQ}{\M{Q}}
\newcommand{\MRr}{\M{R}}

\newcommand{\MT}{\M{T}}
\newcommand{\MU}{\M{U}}
\newcommand{\MV}{\M{V}}
\newcommand{\MW}{\M{W}}

\newcommand{\V}[1]{ \mathbf{#1} }    

\newcommand{\Vo}{\V{0}} 

\newcommand{\Vb}{\V{b}}

\newcommand{\Ve}{\V{e}}

\newcommand{\Vg}{\V{g}}

\newcommand{\Vq}{\V{q}}
\newcommand{\Vr}{\V{r}}
\newcommand{\Vs}{\V{s}}
\newcommand{\Vt}{\V{t}}
\newcommand{\Vu}{\V{u}}
\newcommand{\Vv}{\V{v}}
\newcommand{\Vw}{\V{w}}
\newcommand{\Vx}{\V{x}}
\newcommand{\Vy}{\V{y}}
\newcommand{\Vz}{\V{z}}

\newcommand{\Vgt}{\boldsymbol{\t}}

\newcommand{\Cn}[1]{\mathbb{C}^#1}
\newcommand{\Cmn}[2]{\mathbb{C}^{#1 \times #2}}

\def\eqv{\Leftrightarrow}

\def\rite{\rightarrow}

\def\nm#1{\left\|#1\right\|}
\def\2nm#1{\|#1\|_2}

\def\Ra#1{\mathrm{range}(#1)}

\newcommand{\ars}[1]{\left[ \begin{array}{#1}}
\newcommand{\are}{\end{array} \right] }
\newcommand{\oars}[1]{\begin{array}{#1}}
\newcommand{\oare}{\end{array}}
\newcommand{\rars}[1]{\left( \begin{array}{#1}}
\newcommand{\rare}{\end{array} \right) }

\newcommand{\eqs}{\begin{eqnarray}}
\newcommand{\eqe}{\end{eqnarray}}
\newcommand{\eqsn}{\begin{eqnarray*}}
\newcommand{\eqen}{\end{eqnarray*}}

\newcommand{\bmp}[2]{\begin{minipage}#1{#2}}
\newcommand{\emp}{\end{minipage}}

\newcommand{\ens}{\begin{enumerate}}
\newcommand{\ene}{\end{enumerate}}

\newcommand{\its}{\begin{itemize}}
\newcommand{\ite}{\end{itemize}}

\newcommand{\des}{\begin{description}}
\newcommand{\dee}{\end{description}}

\def\brem{\begin{remark}}
\def\erem{\end{remark}}
\def\defs{\begin{definition}}
\def\defe{\end{definition}}
\def\teos{\begin{theorem}}
\def\teoe{\end{theorem}}
\def\prfs{\begin{proof}}
\def\prfe{\end{proof}}
\def\exas{\begin{exampl}}
\def\exae{\end{exampl}}
\def\excs{\begin{exercise}}
\def\exce{\end{exercise}}
\def\cors{\begin{corollary}}
\def\core{\end{corollary}}

\providecommand{\prn}[1]{\left(#1\right)}

\def\curl#1{\left\{#1\right\}}

\def\ul{\underline}
\def\wh{\widehat}

\def\wtl{\widetilde}
\def\la{\langle}
\def\ra{\rangle}

\newcommand{\comment}[1]{} 

\newcommand{\bea}{\left[ \begin{array} }
\newcommand{\eea}{ \end{array} \right] }

\newcommand{\mrg}{\mbox{Range}}

\usepackage{amsopn}

\begin{document}

\title{A survey of subspace recycling iterative methods \protect\thanks{Document compiled on \today }}

\author[1]{Kirk M. Soodhalter*}

\author[2]{Eric de Sturler}

\author[3]{Misha Kilmer}

\authormark{SOODHALTER \textsc{et al}}

\address[1]{\orgdiv{School of Mathematics}, \orgname{Trinity College Dublin}, \orgaddress{College Green, \state{Dublin}, \country{Ireland}}, ksoodha@maths.tcd.ie}

\address[2]{\orgdiv{Department of Mathematics}, \orgname{Virginia Tech}, \orgaddress{Blacksburg, \state{Virginia}, \country{United States of America}}, sturler@vt.edu}

\address[3]{\orgdiv{Department of Mathematics}, \orgname{Tufts University}, \orgaddress{Medford, \state{Massachussetts}, \country{United States of America}}, misha.kilmer@tufts.edu}

\corres{*School of Mathematics, Trinity College Dublin, College Green, Dublin 2, Ireland. \email{ksoodha@maths.tcd.ie}}


\abstract{
This survey concerns \emph{subspace recycling methods}, a popular class of iterative methods that enable effective reuse of subspace information in order to speed up convergence and find good initial vectors over a sequence of linear systems with
slowly changing coefficient matrices, multiple right-hand sides, or both.
The subspace information that is recycled is usually generated during
the run of an iterative method (usually a Krylov subspace method) on one or more of the systems.
Following introduction of definitions and notation, we examine the history of early augmentation schemes along with deflation preconditioning schemes and
their influence on the development of recycling methods. We then discuss a general residual constraint framework through which many augmented
Krylov and recycling methods can both be viewed. We review several augmented and recycling methods within this framework.
We then discuss some known effective strategies for choosing subspaces to recycle before
taking the reader through more recent developments that have generalized recycling for (sequences of) shifted linear systems, some of them with multiple right-hand sides in mind.
We round out our survey with a brief review of application
areas that have seen benefit from subspace recycling methods.
}


\maketitle

\footnotetext{\textbf{Abbreviations:} CG -- Conjugate Gradients; GMRES -- Generalized Minimum Residual Method; GKB -- Golub-Kahan Bidiagonalization; LSQR -- GKB-based least-squares iterative solver for tall rectangular problems; HPD -- Hermitian positive-definite}

\section{Introduction}\label{sec.intro}

In many applications in the computational sciences, there is a need to solve many hundreds or thousands of large-scale linear systems of the form
\begin{equation}\label{eqn.Aixibi}
	\vek A^{(i)}\vek x_{\ell}^{(i)} = \vek b_{\ell}^{(i)}\qquad i=1,2,\ldots\qquad\ell=1,2,\ldots
\end{equation}
where the systems indexed by $i$ are available in sequence rather than simultaneously, and for each $i$ all right-hand sides indexed by $\ell$
are available simultaneously.
We consider the case where consecutive
coefficient matrices are sufficiently ``closely related'' to exploit the relationship between them.
Such sequences of problems arise in a diverse array of applications; see \Cref{section.applications}.
The coefficient matrices are often sparse or otherwise allow
efficient matrix-vector multiplication
in a matrix-free fashion.
Generally, the dimension
of the matrix is so large that matrix-free iterative methods (e.g., Krylov
subspace methods or multigrid) are the most viable choice for these problems.  With such an iterative method in hand, the most straightforward approach would be
to apply it to each consecutive linear system, with no consideration of the relationship between systems.  However, speedups in convergence
and good initial guesses can be achieved by
exploiting the closeness of consecutive coefficient matrices.

This survey concerns \emph{(Krylov) subspace recycling methods}, a popular class of iterative methods enabling effective reuse of subspace information generated during
the run of an iterative method (usually a Krylov subspace method) applied to $\vek A^{(i)}\vek x^{(i)} = \vek b^{(i)}$, for reuse either after a cycle of iterations (for the same system) or during the iteration applied to
$\vek A^{(i+1)}\vek x^{(i+1)} = \vek b^{(i+1)}$.

\paragraph{Organization of the survey}
In the next section, we briefly review
some basic information about Krylov subspaces methods.
We include some important concepts related to residual projection methods that assist in understanding
recycling strategies.
In \Cref{sec.history-schemes}, we discuss a number of precursors and related techniques.  Some of these are direct forbearers to the current recycling methods, while others
were proposed with different theoretical/practical concerns in mind, but they can be interpreted in the same mathematical framework.
In \Cref{sec.history-recyc}, we give an overview of the state-of-the-art augmentation-based recycling methods.  We describe a general framework,
extending those proposed in \cite{GGL.2013,Gaul.2014-phd,Gutknecht.biCG-aug.2014}, that can be used to described the majority of recycling approaches.
In \Cref{sec.practical-realizations}, we describe generic examples of methods in this general framework for linear systems with
coefficient matrices of various structure and for different residual constraints, and we discuss effective strategies for choosing a subspace to recycle
in \Cref{section.whatrecyc}.
In \Cref{section.additional-structure}, we discuss strategies
to take advantage of recycling for
families of systems with additional structure, in particular, solving multiple shifted systems for each coefficient matrix $\vek A^{(i)}$.  We then briefly discuss in \Cref{section.applications} a variety of scientific, computational,
and engineering applications which have benefited from incorporating a recycling strategy into their solvers.  To conclude, we discuss some challenges which remain.

\paragraph{Notation}  Boldface capital letters denote matrices.  Boldface lowercase letters denote vectors.  Non-boldface letters (Latin and Greek) denote scalar quantities; and, when necessary, we use matching
non-boldface letters to denote entries of a matrix or vector denoted by the same letter (e.g., $h_{ij}\in\bbC$ denotes the entries of the matrix $\vek H$).  The matrices $\vek P$ and $\vek Q$ are used to denote specific projectors that arise in residual projection methods. Calligraphic
letters denote subspaces, and we often use a matching uppercase boldface letter to denote a matrix having that subspace as the span of its columns (e.g., the columns of $\vek U\in\Cmn{n}{k}$ span $\CU\subset\Cn{n}$ which has
dimension $k$).  When these subspaces are without a tilde above (e.g., $\CU$ and $\CV_{j}$), they denote correction spaces from which updates to
solution approximations are drawn.  When such subspaces are written with tildes above (e.g., $\widetilde{\CU}$ and $\widetilde{\CV}_{j}$),
they denote residual constraint spaces used to determine which element is drawn from the correction space; i.e., by enforcing that the new residual
should be orthogonal to the constraint spaces, cf. \eqref{eqn.aug-PG}.
All norms are assumed to be the $2$-norm unless otherwise indicated.

\section{Iterative projection methods -- Krylov Subspace Methods}\label{section.kryl-basic}
Krylov subspace iterative methods are a well-known class of methods for computing an approximation $\vek t_{j}$ to the initial error, $\vek t$, such that for an initial guess, $\vek x_0$, the solution
satisfies $\vek x = \vek x_{0} + \vek t \approx \vek x_0 + \vek t_{j}$; i.e., we are approximating $\vek t$ which solves
\begin{equation}\label{eqn.Axb}
	\vek A\prn{\vek x_{0} + \vek t} = \vek b,\qquad\vek A\in\Cmn{n}{n},\qquad\vek b\in\Cn{n}.
\end{equation}
For solving a linear system of the form \eqref{eqn.Axb} with $\vek r_{0} = \vek b-\vek A\vek x_{0}$,
one builds the Krylov subspace
\begin{equation}\nn
	\CK_{j}(\vek A,\vek r_{0}) = \curl{\vek r_{0},\vek A\vek r_{0},\vek A^{2}\vek r_{0},\ldots,\vek A^{j-1}\vek r_{0}}
\end{equation}
iteratively (at the cost of one matrix-vector product per iteration).
At iteration $j$,
$\vek t_{j}\in\CK_{j}(\vek A,\vek r_{0})$ is selected according to some constraint on the residual $\vek r_{j}$ and
$\vek x_{j} = \vek x_{0} + \vek t_{j}$ is the $j$th approximation.  We call $\Vt_{j}$ the \emph{correction} and the space
from which it is drawn the \emph{correction space}.

\paragraph{An example -- the Generalized Minimum Residual (GMRES) Method}
To illustrate how these methods work, we focus briefly on GMRES \cite{Saad.GMRES.1986},
in which we select
\begin{equation}\label{eqn.pg-GMRES}
\vek t_{j}\in\CK_{j}(\vek A,\vek r_{0})\msuchthat \vek r_{j}\perp\vek A\CK_{j}(\vek A,\vek r_{0}).
\end{equation}
Methods with such a residual orthogonality constraint are called \emph{residual projection methods}, because
the constraint defines the projection (oblique or orthogonal) of the residual onto
$\CK_{j+1}(\MA,\Vr_0)$, and we call $\vek A\CK_{j}(\vek A,\vek r_{0})$ here the \emph{constraint space}.
Characterization of these methods according
to residual constraints is helpful in understanding the general structure of recycling methods.
This particular constraint is equivalent to solving the residual minimization problem
\begin{equation}\label{eqn.gmres-full-min}
	\vek t_{j} = \argmin{\boldsymbol\tau\in\CK_{j}(\vek A,\vek r_{0})}\nm{\vek b - \vek A(\vek x_{0} + \boldsymbol\tau)} ,
\end{equation}
and leads to the approximation
\begin{equation}\label{eqn.GMRES-resid-proj}
	\vek t_{j} = \vek P_{\CK_{j}}\vek t\mand \vek r_{j} = \prn{\vek I - \vek Q_{\CK_{j}}}\vek r_{0},
\end{equation}
where $\vek P_{\CK_{j}}$ is the $\prn{\vek A^{\ast}\vek A}$-orthogonal projector onto $\CK_{j}\prn{\vek A,\vek r_{0}}$, and
$ \vek Q_{\CK_{j}}$ is the orthogonal projector onto $\vek A\CK_{j}\prn{\vek A,\vek r_{0}}$.

\paragraph{The Arnoldi orthogonalization procedure leads to a practical implementation}
The Arnoldi process builds an orthonormal basis for $\CK_{j+1}(\vek A,\vek r_{0})$. Set $\Vv_1 = \b^{-1}\Vr_0$ with $\b = \|\Vr_0\|$. At iteration $j$, we compute
\eqs \label{eqn.Arnoldi.process}
  \Vv_{j+1}h_{j+1,j} & = & \MA\Vv_j - \sum_{i=1}^{j} \Vv_i h_{i,j},
  \quad \mbox{with} \quad h_{i,j} = \Vv_{i}^{*} \MA \Vv_j \mbox{ and }
  h_{j+1,j} = \| \MA\Vv_j - \sum_{i=1}^{j} \Vv_i h_{i,j} \| ,
\eqe
so that $\Vv_{j+1}$ is a unit vector. The process computes the matrices
\begin{equation}\label{eqn.arnoldi-vectors}
	\vek V_{j+1}=\begin{bmatrix} \vek v_{1} & \vek v_{2} & \cdots & \vek v_{j+1} \end{bmatrix}\in\Cmn{n}{(j+1)}\mand \underline{\vek H}_{j}\in\Cmn{(j+1)}{j},
\end{equation}
such that $\MV_{j+1}^*\MV_{j+1} = \MI$, $\Ra{\MV_{j+1}}=\CK_{j+1}(\vek A,\vek r_{0})$,
$\underline{\vek H}_{j}$ is an upper Hessenberg matrix with components $h_{i,j}$, and
\begin{equation}\label{eqn.arnoldi-relation}
	\vek A\vek V_{j} = \vek V_{j+1}\underline{\vek H}_{j} = \vek V_{j}\vek H_{j} + h_{j+1,j}\vek v_{j+1}\vek e_{j}^{*},
\end{equation}
where $\vek H_{j}\in\Cmn{j}{j}$ is simply the first $j$ rows of $\underline{\vek H}_{j}$.
Using \eqref{eqn.arnoldi-relation}, one can reduce the minimization \eqref{eqn.gmres-full-min} to a smaller $(j+1)\times j$
least-squares minimization problem
\begin{equation}\label{eqn.gmres-ls-min}
	\vek y_{j} = \argmin{\vek y\in\bbC^{j}}\nm{\underline{\vek H}_{j}\vek y - \beta\vek e_{1}}\mand \vek t_{j} = \vek V_{j}\vek y_{j}.
\end{equation}

\paragraph{The Hermitian Lanczos process}
In this case, we can use the Arnoldi process above,
while exploiting that $\vek A$ is Hermitian and hence
$\MH_j = \vek V_j^* \vek A \vek V_j$ is itself Hermitian. It follows that $\MH_j$ is tridiagonal and $\vek V_j$ can be computed efficiently with the three-term recurrence,
\begin{equation}\label{eqn.Lanc_rec-1}
  \vek v_{j+1}h_{j+1,j} = \vek A\vek v_{j} - \vek v_{j}h_{j,j} - \vek v_{j-1} h_{j-1,j},
      \qquad \mbox{ and }\qquad \vek v_2 h_{2,1} = \vek A\vek v_1 - \vek v_1 h_{1,1}
 \end{equation}
with the $h_{i,j}$ as in (\ref{eqn.Arnoldi.process})
and $h_{j-1,j} = h_{j,j-1}$ available from the previous iteration.
This recurrence is called the Lanczos recurrence, leading to the Lanczos relation (with
$\vek T$ for tridiagonal)
\eqs \label{eqn.Lanc_rec-2}
  \vek A \vek V_j & = & \vek V_{j+1}\underline{\vek T}_j =
  \vek V_j \vek T_j + \vek v_{j+1} h_{j+1,j} \Ve_{j}^*.
\eqe
The short recurrence \eqref{eqn.Lanc_rec-1}
leads to a great reduction in the memory
requirements as
we only need to store the two most recently generated Lanczos vectors.
Using the thin QR-decomposition
$\ul{\MT}_j = \ul{\MG}_j^{(j+1) \times j} \MRr_j$ to solve
\eqref{eqn.gmres-ls-min}, we can write the solution update at
step $j$ as $\Vt_j = \MV_j \MRr_j^{-1} \ul{\MG}_j^{*} \Ve_1 \beta$.
Since the coordinate vector $\wtl{\Vy} = \ul{\MG}_j^{*} \Ve_1 \beta$
changes only in its last coefficient from one iteration to the next, this leads to an
efficient update procedure as follows.
Performing the \emph{change of basis} $\MW_j = \MV_j \MRr_j^{-1}$ leads to the additional 3-term
recurrence $\Vw_j = r_{j,j}^{-1}(\Vv_j - \Vw_{j-1}r_{j-1,j} - \Vw_{j-2}r_{j-2,j})$ and the update
$\Vx_j = \Vx_0 + \Vt_j = \Vx_0 + \MW_j \wtl{\Vy} = \Vx_{j-1} + \Vw_j \wtl{y}_j$.
The resulting method is called MINRES \cite{Paige1975}.

If we assume additionally that $\vek A$ is Hermitian positive definite (HPD),
we can minimize the error in the $\MA$-norm, $\| \vek t-\vek t_{j}\|_{\vek A}$, which corresponds to enforcing the residual constraint $\vek r_{j}\perp \CK_{j}\prn{\vek A,\vek r_{0}}$ \cite[Sections 6.4 and 6.7.1]{Saad.Iter.Meth.Sparse.2003}.
Considering again the Lanczos relation (\ref{eqn.Lanc_rec-2}), we have that
$\MT_j$ is also HPD and allows the LU decomposition
$\MT_j = \ML_j \MU_j$ with unit bidiagonal $\ML_j$.
Now taking the change of basis $\MW_j = \MV_j \MU_j^{-1}$ and
$\wtl{\Vy} = \ML_j^{-1} \Ve_1$ and eliminating the explicit
Lanczos recurrence using that $\Vv_{j+1}$
is just a normalization of the residual $\Vr_{j}$, we 
obtain the celebrated method of conjugate gradients \cite{Hestenes.Stiefel.CG.1952}.

\paragraph{The non-Hermitian Lanczos process}
For general non-Hermitian linear systems, there exist short recurrence methods, though they generally do not lead to
an orthonormal basis for the Krylov subspace \cite{Faber.Manteuffel.1984}.
Instead, one simultaneously generates dual bases
for the subspaces $\CK_{j} (\vek A,\vek r_{0})$ and $\CK_{j}(\vek A^{\ast},\wh{\vek r}_{0})$ where
$\wh{\vek r}_{0}$ is either the initial residual of a dual problem involving $\vek A^{\ast}$, $\vek r_{0}$ itself, or some other
non-zero vector. The biorthogonal Lanczos process is a short recurrence for iteratively generating these bases simultaneously
at the cost of
one application of $\vek A$ and one of $\vek A^{\ast}$ per iteration,
\begin{equation}\label{eqn.biortho.bases}
{\rm span}\left\lbrace \vek v_{1}, \vek v_{2}, \ldots, \vek v_{j} \right\rbrace = \CK_{j}(\vek A,\vek r_{0}) \qquad\mbox{and}\qquad {\rm span}\left\lbrace \wh{\vek v}_{1}, \wh{\vek v}_{2}, \ldots, \wh{\vek v}_{j} \right\rbrace = \CK_{j}(\vek A^{\ast},\wh{\vek r}_{0}).
\end{equation}
The basis vectors are constructed to satisfy $\vek v_{\ell}^{\ast}\wh{\vek v}_{k}=\delta_{\ell k}$ (hence they are \emph{biorthogonal}).
The biorthogonal Lanczos relations are as follows, where $\MT_j, \wh{\MT}_j$ are tridiagonal:
\begin{equation}\label{eqn.biorth-lanzcos}
	\vek A\vek V_{j} = \vek V_{j}\vek T_{j}+h_{j+1,j}\vek v_{j+1}\vek e_{j}^{\ast},\qquad\mbox{and}\qquad \vek A^{\ast}\wh{\vek V}_{j} =
\wh{\vek V}_{j} \wh{\MT}_{j} +\wh{h}_{j+1,j}\wh{\vek v}_{j+1}\vek e_{j}^{\ast} .
\end{equation}
The Petrov-Galerkin condition $\vek r_{j}\perp \CK_{j}(\vek A^{\ast},\widetilde{\vek r}_{0})$, giving rise to the biconjugate gradient (BiCG) method \cite{Fletcher-BiCG.1976,Lanczos.nonsym.1952ROHTUA}, requires
\begin{equation}\label{eqn.bicg}
	\mbox{Solve}\qquad \vek T_{j}\vek y_{j} = \beta\vek e_{1}\qquad\mbox{and set}\qquad\vek t_{j} = \vek W_{j}\vek y_{j}.
\end{equation}
Similar to \eqref{eqn.arnoldi-relation} for Hermitian systems,
\eqref{eqn.biorth-lanzcos} allow us to define efficient short recurrence iterative methods for non-Hermitian systems using the LU decomposition
$\MT_j = \ML_j \MU_j$, change of basis $\MW_j = \MV_j\MU_j^{-1}$,
and setting $\wtl{\Vy} = \ML_j^{-1} \prn{\beta \Ve_1}$
(in general, $\MT_j$ and $\wh{\MT}_j$ are closely related, and only one LU-decomposition is computed).
An analog to GMRES/MINRES also exists in this setting, called the quasi minimum residual (QMR) method \cite{Freund.QMR}.
At iteration $j$, this is leads to
\begin{equation}\label{eqn.qmr}
	\mbox{Minimize } \left\|\ul{\vek T}_{j}\vek y_{j} - \beta \vek e_{1}\right\|,\qquad\mbox{and set}\qquad\vek t_{j}=\vek W_{j}\vek y_{j}.
\end{equation}
In cases where the action of $\vek A^{\ast}$ is unavailable, so-called \emph{transpose free} variants are available; see, e.g.,
\cite{Sonneveld.CGS.1989,Freund.TFQMR,vandervorst.bicgstab.1992}.
The BiCGStab($\ell$) \cite{bicgstabl.1993} method was introduced to stabilize the oscillatory convergence pattern often exhibited by
BiCG \cite{vandervorst.bicgstab.1992}.  This is accomplished by alternating between $\ell$ steps of BiCG and an
$\ell$-cycle of GMRES, effectively building a hybrid residual polynomial from the BiCG and GMRES polynomials.
Another related method is IDR($s$) \cite{Sonneveld2008}, which was shown to fit into
the Petrov-Galerkin residual projection framework \cite{Simoncini2010a} and also to be a generalization of BiCGStab($\ell$).

\section{Subspace augmentation and the strategy of deflation}\label{sec.history-schemes}
Subspace recycling extends ideas from the last few decades for the preservation of
information between iteration cycles or between different linear systems.
This is done to mitigate the effects of discarding basis vectors due to memory requirements
as well as to accelerate the convergence of an iterative method.
A number of acceleration strategies occupy the same or a similar theoretical framework as augmentation-based subspace recycling, and
we briefly touch upon these as well.

\textbf{What we mean by augmentation}: In this survey, we use the terms \emph{subspace augmentation} and \emph{augmented method} to refer to
any iterative method (here, a Krylov subspace method) which uses a sum correction space $\CU + \CV_{j}$, where $\CV_{j}$ is the space generated
by the iterative method, and $\CU$ is the fixed \emph{augmentation space}.  Later, we discuss cases where $\CV_{j}$ is a Krylov subspace, and one can consider
the case in which the Krylov subspace is generated by the coefficient matrix/residual pair $\prn{\vek A,\vek r_{0}}$
versus cases in which the Krylov subspace is generated by a projected pair of the form
$\prn{\prn{\vek I - \vek Q}\vek A,\prn{\vek I - \vek Q}\vek r_{0}}$, cf. \Cref{sec.history-recyc}.

\paragraph{Warm restarting and simple information reuse}
The simplest information-reuse approach is called \emph{warm restarting}.
As we generate approximate solutions for the sequence of problems
(\ref{eqn.Aixibi})
we use the approximate solution
for $\vek A^{(i)}\vek x^{(i)} = \vek b^{(i)}$ as the initial approximation $\vek x_{0}$ for
the system $\vek A^{(i+1)}\vek x^{(i+1)} = \vek b^{(i+1)}$.  As consecutive systems in (\ref{eqn.Aixibi}) are
close (e.g., consecutive systems in a nonlinear optimization scheme), the 
approximate solution to system $i$ is likely a high-quality
approximation for the solution to the system $i+1$.

For nonsymmetric systems, the first mention of a (non-trivial) strategy to select direction vectors for orthogonalization
(to the best of our knowledge) is in \cite{JR.1985}, where
the authors performed a numerical study of algorithms associated with preconditioned conjugate gradients.
The authors suggest that for non-symmetric problems and a fixed (but sufficient) length recurrence, convergence is improved by retaining the search directions
associated with the (relative) largest orthogonalization coefficients.

\paragraph{Augmentation via flexible preconditioning}
Augmentation of a Krylov subspace was proposed in a 1997 paper, which presented the idea in the framework of the flexible GMRES
method  \cite{Chapman1997}.
Flexible GMRES is a modification of right-preconditioned GMRES that accommodates the use of a different preconditioner at each step.
Recall that for a fixed right preconditioner $\vek M$, one applies
an iterative solver to the problem $\vek A\vek M^{-1}\vek y = \vek b$ where $\vek y = \vek M\vek x$,
and for an approximation $\widetilde{\vek y}$ one generates the approximation 
$\widetilde{\vek x} = \vek M^{-1}\widetilde{\vek y}$.

Flexible preconditioning complicates this situation.
In flexible GMRES, one must store the Arnoldi vectors from (\ref{eqn.arnoldi-vectors})
for a Krylov subspace  $\CK_{j}\prn{\vek A\boldsymbol{\mathcal{M}}^{-1},\vek r_{0}}$, where $\boldsymbol{\mathcal{M}}$ is an unknown implicitly induced
preconditioner\footnote{It was pointed out in \cite{Notay.2000}
that flexible preconditioning implicitly induces exact an preconditioner
$\boldsymbol{\mathcal{M}}$ defined
by $\boldsymbol{\mathcal{M}}^{-1}\vek v_{i} = \vek M_{i}^{-1}\vek v_{i}$, for $i=1,2,\ldots, j$.
}, as well as a set of flexibly preconditioned Arnoldi vectors
$\vek Z_{j} = \begin{bmatrix}\vek z_{1} & \vek z_{2 }& \ldots &\vek z_{j} \end{bmatrix}$, with
$\vek z_{i} = \vek M_{i}^{-1}\vek v_{i}$, spanning $\boldsymbol{\mathcal{M}}^{-1}\CK_{j}\prn{\vek A\boldsymbol{\mathcal{M}}^{-1},\vek r_{0}}$.
One solves the GMRES minimization (\ref{eqn.gmres-ls-min}) for $\vek y_{j}$
as before but then updates $\vek x_{j} = \vek x_{0} + \vek Z_{j}\vek y_{j}$.

The authors of \cite{Chapman1997} pointed out that one can use this
framework to augment an existing Krylov subspace.
Assume we have run $j-k$ iterations of GMRES and
$\vek U = \begin{bmatrix}\vek u_{1} & \vek u_{2} & \cdots & \vek u_{k}\end{bmatrix}\in\Cmn{n}{k}$
spans a $k$ dimensional subspace $\CU$
to use in concert with the Krylov subspace for a correction. For the next $k$ iterations, one can
implicitly define the action of unspecified flexible preconditioners $\vek M_{j-k+i}$ with the mappings $\vek M_{j-k+i}^{-1}\vek v_{j-k+i} = \vek u_{i}$ for
$i=1,2,\ldots, k$.
This defines a version of flexible GMRES that minimizes the residual using the constraint space
$\vek A\CK_{j-k}\prn{\vek A,\vek r_{0}} + \vek A\,\CU$.  This minimization is \textbf{not} over a direct sum of spaces.  The lack of orthogonality
between the subspaces 
reduces its effectiveness.  A similar performance penalty due to lack of orthogonality
was the impetus for the development of the GCROT method \cite{deSturler.GCROT.1999}.  Augmentation-type methods for ill-posed problems
have been developed that rest on similar ideas \cite{DGH.2014}.

\paragraph{Deflation}
Deflation type methods use some of the same underlying mathematical principles as recycling methods, but they are couched in a slightly different language and aimed at different goals.
Deflation techniques were originally proposed to treat a single system via left preconditioning using a projector onto the orthogonal
complement of some subspace \cite{Nico_87}. The convergence is then determined by a projected operator and right-hand side, and projecting away from an appropriate subspace (e.g., some invariant subspace of the coefficient matrix),
leads to faster convergence.

In \cite{Erl.Nabben.Defl-balPrec.2008}
the authors present a technique to analyze methods such as the deflated restarting method of Morgan \cite{Morgan.GMRESDR.2002} (discussed in greater detail in
\Cref{sec.history-recyc}).  They interpret such deflated methods (specifically balancing preconditioners \cite{Mandel.BalPrec.1993})
as being preconditioned with \emph{deflation preconditioners} that project into an associated augmentation space.
For example, suppose we are solving
(\ref{eqn.Axb}), where $\vek t$ represents the initial error.  Let $\MU \in \Cmn{n}{k}$ and
$\CU = \Ra{\MU}$ be the subspace of dimension $k$ from which we want to construct
an approximation to $\vek t$. In addition, let $\wtl{\MU} \in \Cmn{n}{k}$ and
$\wtl{\CU} = \Ra{\wtl{\MU}}$ be the constraint space, so that $\Vt = \MU \Vs$
satisfies $\Vr = \Vb - \MA(\Vx_0 + \MU\Vs) \perp \wtl{\CU}$ or equivalently
$\wtl{\MU}^*\Vr = \wtl{\MU}^*(\Vb - \MA(\Vx_0 + \MU\Vs)) = \Vo$.
$\vek P$ denotes the corresponding (oblique or orthogonal) projector onto $\CU$,
that is, $\MP \Vt = \MU \Vgt$, which implies
$\vek P = \vek U\prn{\widetilde{\vek U}^{\ast}\vek A\vek U}^{-1}\widetilde{\vek U}^{\ast}\vek A$,
and $\vek Q$ denotes the sibling projector\footnote{
We call these sibling projectors because specifying one of them in a projection method determines the other unambiguously. Defining a projection method via a specific projector $\vek P$ being applied to the error also determines $\vek Q$.  Vice-versa, defining a method
via the application of a projector $\vek Q$ to the residual determines $\vek P$.
}
onto $\vek A\CU$ that satisfies
\begin{equation}\label{eqn.proj-relation}
	\vek Q\vek A = \vek A\vek P.
\end{equation}
The deflation-preconditioning process splits the initial error:
\begin{equation}\label{eqn.error-splitting}
\vek t = \vek P\vek t + \prn{\vek I - \vek P}\vek t .
\end{equation}
Since $\vek A\vek t = \vek r_{0}$, we can directly compute $\vek P\vek t = \vek U\prn{\widetilde{\vek U}^{\ast}\vek A\vek U}^{-1}\widetilde{\vek U}^{\ast}\vek r_{0}$, which can be interpreted as an approximation of the initial error in the space $\CU$.
It remains then to approximate $\prn{\vek I - \vek P}\vek t$, for which we have the following system
of equations,
\eqs
\nonumber
  \MA (\Vx_0 + \Vt) = \MA (\Vx_0 + \MP\Vt + (\MI - \MP)\Vt ) = \Vb & \eqv &
  \MA(\MI - \MP)\Vt = \Vb - \MA \Vx_0 - \MA \MP \Vt = \Vr_0 - \MA\MP\Vt \quad \eqv \\
\label{eqn.derive-proj-Axb}
  \MA(\MI - \MP)\Vt = \Vr_0 - \MQ\MA\Vt = (\MI - \MQ)\Vr_0 & \eqv &
  (\MI - \MQ)\MA\Vt = (\MI - \MQ)\Vr_0 .
\eqe
Hence, we solve the deflation-preconditioned linear system
\begin{equation}\label{eqn.proj-Axb}
	\prn{\vek I - \vek Q}\vek A\vek t = \prn{\vek I - \vek Q}\vek r_{0}.
\end{equation}
If we apply an iterative method to (\ref{eqn.proj-Axb}), which returns an approximation $\vek t_{j}$ at iteration $j$, then the full approximation at this step is $\vek x_{j} = \vek x_{0} + \vek P\vek t + \prn{\vek I - \vek P}\vek t_{j}$.
The residual of this approximation satisfies $\vek r_{j} = \vek b - \vek A\vek x_{j} = \prn{\vek I - \vek Q}\prn{\vek b - \vek A\vek x_{j} }$. Thus it is equal to the residual of the projected problem,
meaning residual convergence is completely determined by properties of the
projected problem (\ref{eqn.proj-Axb}).  In the Ph.D. thesis \cite{Gaul.2014-phd}, it was confirmed that GCRO-based recycling schemes are indeed equivalent to particular deflation schemes for certain choices of \emph{orthogonal} projectors $\vek P$ and $\vek Q$.

Methods such as that in
\cite{EigCG-nonsymm-lanczos,MR2599785} approach the same strategy by performing a one-time projection away from specific accurately-computed
eigenvectors to remove their influence from the iteration.  This does not quite fit into the recycling framework as we have described it in this survey;
but, as with deflation, it is based on some of the same ideas.
\paragraph{Multigrid and domain decomposition}

The correction $ \vek P\vek t + \prn{\vek I - \vek P}\vek t_{j}$ comes from the space
$\CU + \CK_{j}\prn{\prn{\vek I - \vek Q}\vek A,\prn{\vek I - \vek Q}\vek r_{0}}$.
This augmented subspace interpretation of deflation-based preconditioning can be widened to also include methods involving restrictions onto and prolongation from subgrids or coarser grids such as domain decomposition and multigrid.
This has been observed, e.g., in \cite{dolean2015introduction}.
Consider, for example, a simple restriction-solution-prolongation procedure, such as a one-level multigrid V-cycle or a domain decomposition subdomain solve.   Let $\mathfrak{R}\in\Cmn{k}{n}$ be the restriction operator associated to the procedure, and
let $\mathfrak{P}\in\Cmn{n}{k}$ be the associated prolongation operator.  Then the restriction-solution-prolongation procedure applied to (\ref{eqn.Axb}) can be represented as
\begin{equation}\label{eqn.restrict-prolong}
	\widehat{\vek x} = \vek x_{0} + \mathfrak{P}\prn{\mathfrak{R}\vek A\mathfrak{P}}^{-1}\mathfrak{R}\vek r_{0} \qquad{\rm and}\qquad \widehat{\vek r} = \vek r_{0} - \vek A\mathfrak{P}\prn{\mathfrak{R}\vek A\mathfrak{P}}^{-1}\mathfrak{R}\vek r_{0}.
\end{equation}
If we further assume that $\vek A$ is Hermitian positive-definite, and that the restriction and prolongation operators satisfy the relationship $\mathfrak{R}=c\mathfrak{P}^{\ast}$, where $c\in\bbC$,
then (\ref{eqn.restrict-prolong}) reduces to
\begin{equation}\nonumber
	\widehat{\vek x}\leftarrow \vek x_{0} + \vek P\vek t \qquad{\rm and}\qquad \widehat{\vek r} \leftarrow\vek r_{0} - \vek Q\vek r_{0}.
\end{equation}
where $\vek P,\,\vek Q$ are projectors of the form $\vek P = \mathfrak{P}\prn{\mathfrak{P}^{\ast}\vek A\mathfrak{P}}^{-1}\mathfrak{P}^{\ast}\vek A$ and $\vek Q = \vek A\mathfrak{P}\prn{\mathfrak{P}^{\ast}\vek A\mathfrak{P}}^{-1}\mathfrak{P}^{\ast}$.  These
do satisfy (\ref{eqn.proj-relation}), and can be interpreted in the deflation-augmentation framework, where the augmentation space is now some sort of interpolation or restricted subdomain space. Indeed, it follows from the fact that this restriction-solution-prolongation minimizes the error
in the $\vek A$ -norm that $\CU = {\rm range}\prn{\mathfrak{P}} = \widetilde{\CU}$.
The projectors $\vek P$ and $\vek Q$ are the $\vek A$-orthogonal projector onto ${\rm range}\prn{\mathfrak{P}}$ and the $\vek A^{-1}$-orthogonal projector onto ${\rm range}\prn{\vek A\mathfrak{P}}$, respectively.  For a general matrix $\vek A$, the procedure becomes non-optimal with $\vek P$ being
the oblique projector onto $\CU$ along $\prn{\vek A^{\ast}\CU}^{\perp}$, and $\vek Q$ is the projector onto $\vek A\CU$ along $\CU^{\perp}$ (see \cite{dolean2015introduction} and references therein).

\paragraph{Bordering methods}
Bordering methods were originally introduced as a means to solve a linear system by augmenting its rows and columns in a particular way to induce a singular consistent
problem whose family of solutions have within them the solution to the original problem \cite{fadeev1963computational,axelsson1996application}.
A straightforward explanation can be found in the presentation \cite{Neytcheva.Deflation.Talk}.  Looking at the details, one sees elements of
something that looks like augmentation and is also a precursor to deflation preconditioning.

\paragraph{Polynomial preconditioning}
Some types of polynomial preconditioning and hybrid approaches
(see e.g., \cite{MR1410697,MR1261463,NachReiTref_1992, SimonGallo_1995}) have a strong connection to deflation methods and augmentation for a single linear system.
In such methods, one runs a cycle of GMRES for a given system.  After this cycle, a fixed polynomial is constructed implicitly using information
from the Krylov subspace.  This polynomial can either be used with a stationary iterative method or as a preconditioner for a iterative method such as GMRES.

\section{Subspace recycling}\label{sec.history-recyc}

\paragraph{When is an augmented method subspace recycling?}
In this survey, we focus on a specific type of augmented Krylov subspace methods called \emph{Krylov subspace recycling methods}, a moniker that communicates that the additional correction subspace was {\em recycled} from a previously generated but discarded subspace for its importance in obtaining a fast convergence rate or good initial vectors.
For example, at the end of a GMRES cycle, before restarting,
retaining a subspace of the generated Krylov space may accelerate convergence in the next cycle.
Or, the recycled subspace may have been determined
to damp the influence of
certain parts of the spectrum of the operator \cite{Morgan.GMRESDR.2002,Parks.GCRODR.2006}.
\emph{An important aspect of subspace recycling is that, in addition to augmentation, the Krylov space is changed to complement the recycle space }.
As with deflation, the new Krylov space is generated by the operator
$(\MI-\MQ)\MA$ and a projected residual or right hand side.
The importance of this was observed in \cite{EDS.GCRO.1996} proposing the GCRO method(s); this generally avoids the (possible) stagnation problems of restarted GMRES, generates the right search
space given the existing correction space, and computes the optimal solution over the sum of the correction recycle space and the generated Krylov space.

The main ideas behind Krylov subspace recycling arose from the problem that optimal Krylov methods for a fixed nonsymmetric system need all Krylov iteration (or direction)
vectors in each iteration. This leads to excessive storage requirements (linear in the number of iterations) and computational work (quadratic in the number of iterations). A practical solution is to restart every $m$ iterations or to truncate the set of iteration vectors and orthogonalize only against the last $m$, previous iterations, leading to methods such as GMRES(m) \cite{Saad.GMRES.1986}, truncated GCR(m) \cite{EisElS83}, DIOM and DQGMRES\cite{Saad.Iter.Meth.Sparse.2003}. While these strategies often work well, they can lead to very slow convergence.

\paragraph{Augmentation based on ideas of deflation}
Based on the ideas by Nicolaides for deflated CG \cite{Nico_87}, several researchers, e.g.,
Morgan \cite{Morgan.Restarted-GMRES-eig.1995,Morgan2000}, Karchenko and Yeremin \cite{KarYere1995}, Erhel and collaborators \cite{ErhelBurraPohl_96,BurrErhe1998}, Baglama et al. \cite{Bagl-etal_1998}, and Frank and Vuik \cite{FranVuik_2001}  proposed to use augmentation
with (approximate) eigenvectors to maintain good convergence after restarting. Vuik and collaborators have also considered deflation preconditioned methods based on a priori considerations for good deflation spaces, based on algorithmic or application information, \cite{Shei-etal_2016,Jons-etal_2012}.
These methods may use augmentation or a preconditioning approach to achieve the desired effect.
Subsequently, Morgan proposed several improvements to the augmentation approach in \cite{Morgan.Restarted-GMRES-eig.1995},
leading to the elegant GMRES-DR method that exploits implicit restarting of the Arnoldi recurrence to use deflation \cite{Morgan.GMRESDR.2002}.
Just as for GCRO, rather than augmenting
the standard Krylov subspace generated by GMRES, as in
\cite{Morgan.Restarted-GMRES-eig.1995},
Morgan proposes to augment a Krylov subspace generated by
$(\MI - \MQ)\MA$, cf.  (\ref{eqn.proj-Axb}).
Similar ideas based on implicit restarting and deflated restarts have been proposed for Lanczos-based methods \cite{MR1978154,MR1270124}
and Golub-Kahan Bidiagonalization-based methods (GKB) such as LSQR \cite{BRR-augLSQR.2013,MR3123847,MR2201173}.  There are many other approaches which build on this type of deflated approach which incorporate
flexible preconditioning \cite{GGPV.FGMRES-DR.2012}, simultaneously treat block Krylov methods 
\cite{Morg_2005} with inexact breakdown \cite{Giraud.BGMRES-DR-Inexact.2014}, and use deflation in a
flexibly preconditioned CG approach \cite{CGLV.FGCRODR.2012}.

\paragraph{Optimal augmentation built on GCR, i.e., GCRO-based methods}
In a parallel development, independently, Eirola and Nevanlinna\cite{Eirola.Nevanlinna.1989} proposed a splitting-based iteration in which the splitting is updated in each step with a specially chosen, rank-one matrix. Although the authors did not state this, with the right choices the method is equivalent to GCR \cite{EisElS83}, as observed in \cite{VV.GMRESr.1994}. The GMRESR method proposed in the latter paper\footnote{The general scheme is referred to as GMRES* in \cite{vanDerVorst-book.2009} and only with GMRES(m) as inner method as GMRESR.} replaces the residual as search direction in GCR with an approximation to the error computed by GMRES (or by another iterative method), leading to an inner-outer iterative method.
In \cite{EDS.GCRO.1996} it was observed that optimality over the direct sum of the inner and outer correction spaces could be obtained by maintaining the right orthogonality relations. The resulting method is
called GCRO. As with GMRESR, the GCRO method can be combined with any inner method, such as BiCGStab \cite{EDS.GCRO.1996}. The augmentation
space for GCRO is spanned by the corrections computed by
the sequence of inner iterations. The GCROT method extends GCRO by computing an optimal subspace to recycle for subsequent iterations. The optimality is based on considering the canonical angles between the subspaces generated by restarted GMRES \cite{deSturler.GCROT.1999}.
An extension of GCRO for a sequence of linear systems with a fixed matrix and multiple right hand sides was presented in \cite{EdS_HH96}. In \cite{Parks.GCRODR.2006}, this was further
extended to a sequence of linear systems where the matrix
changes slowly with right hand sides that may or may not be close, leading to the
recycling GCROT and GCRODR methods. In
\cite{KilmerdeSturler2006}, extensions include
(a) a recycling version of MINRES,
(b) utilizing approximate solutions in the recycle space
to get good initial vectors, and (c) using recycling for
a sequence of matrices with a number of shifts for each matrix.
Additional innovations for multiple shifted systems are discussed in
section \ref{section.additional-structure}.
More efficient versions of recycling MINRES, especially including efficient ways
of computing and updating recycle spaces, were proposed in
\cite{Wang.TopOptRecyc.2007,Recyc-Imp-Tom.2010}.
Extensions for recycling in BiCG \cite{Fletcher-BiCG.1976},
BiCGStab \cite{vandervorst.bicgstab.1992},
and IDR(s) \cite{Sonneveld2008} have also been
proposed
\cite{Ahuja.Recyc-BiCG.2012,Ahuja.Recyc-BiCGStab-MOD.2015,NeuGreif_18}.
The idea to use projection on a search space of old solutions for initial vectors
was also proposed
in \cite{Fischer.init-guess-from-prev.1998}, but this paper does
not involve recycling or augmentation.

\paragraph{Analysis}
There has not been a great deal of convergence analysis of subspace recycling methods. However, several papers illuminate certain aspects of the
behavior of these methods and related augmentation and acceleration strategies.  An early analysis focuses on the non-optimal
augmentation strategy \cite{Saad.Deflated-Aug-Krylov.1997}. An analysis of acceleration strategies based on the principle angles
between subspaces provides a comparison of strategies but
not quantitative convergence bounds can be found in  \cite{Eiermann2000}.
There has been some analysis
of using approximate invariant subspaces to accelerate convergence.
This was touched upon in a presentation at the 2012 Householder Symposium and the associated
abstract \cite{desturler.recycling-convergence-abstract.2012}, and it will be elaborated further in our review paper \cite{dSKS.2018}. Another important
way to address these questions might be to use the analysis of the onset of superlinear convergence in a GMRES iteration \cite{Simoncini2005}.
Another justification supported by analytic results can be found
in \cite{KilmerdeSturler2006}.  See \Cref{section.whatrecyc} for further details.
There is also analysis
from the point of view of deflation preconditioning (e.g.,\cite{Erl.Nabben.Defl-balPrec.2008} and \cite{Gaul.2014-phd}), and there has been a related analysis
of the use of approximate deflation preconditioners constructed using previously-generated Arnoldi vectors \cite{Sifuentes2011}.

\paragraph{Alternative approaches}
Independently of augmentation-based subspace recycling, recycling of information was the primary point of \emph{seed-based methods}.
These are not deflation methods.  Rather, they use the whole Krylov space generated for one system to then solve multiple simultaneously available systems.
The block seed methods in \cite{Saad-seed,CN.GalProjAnalMRHS.1999,SimonGallo_1995, Chan1997,KMR.2001} use a block Krylov subspace to solve a subset of the systems and update the remaining
systems using the just-generated block Krylov subspace according to some projection or minimization.  It should be noted, though, that one strategy pursued for combining subspace recycling with strategies for solving
multiple shifted linear systems (cf. \Cref{sec:success.recyc.shift}) built upon these methods and the analysis thereof \cite{S.2014.2}.

For groups of HPD matrices that are all close to one another, a method was proposed which is close in spirit to a recycled CG method \cite{Risler1998}.
The authors propose to reuse the entire Krylov space generated for one coefficient matrix as the augmentation space for the next, which is computationally quite expensive, and the authors suggest to restart once memory has been
exhausted. The method was extended with the use of converged dominant Ritz vectors (related to the largest eigenvalues) \cite{Rey1998} resulting in modest improvements. Both methods, as proposed, require an expensive full recursion in spite of the operator being HPD.
Later, alternative subspaces for augmentation using the theory in \cite{SluiVors_86} were proposed \cite{Gosselet.etal.reuse-Kryl-nonlin.2012}, following the recycling approach.

\subsection{A framework}\label{section.framework}
We now briefly present a general residual constraint framework through which many augmented Krylov subspace
methods can be viewed.  For a more complete view of this framework, see \cite{dSKS.2018} which builds an understanding of these methods
in terms of residual constraints on top of the existing work in \cite{Gaul.2014-phd, GGL.2013, Gutknecht.biCG-aug.2014}.  Consider two subspaces
$\CU,\,\CV_{j}\subset\Cn{n}$, where $\CU$ has fixed dimension $k$ and $\CV_{j}$ is generated by an iterative process such that at step $j$,
$\dim\CV_{j}=j$.  We take $\CU+\CV_{j}$ as our correction space.  We can similarly construct a constraint space  $\widetilde{\CU}+\widetilde{\CV}_{j}$, with
$\widetilde{\CU},\,\widetilde{\CV}_{j}\subset\Cn{n}$ being of fixed dimension $k$ and iteratively generated with dimension $j$ at iteration $j$, respectively.
We assume for simplicity that
$\dim\prn{\CU+\CV_{j}} = \dim\prn{\widetilde{\CU}+\widetilde{\CV}_{j}} = k+j$ (direct sums).  Then the general augmented projection method becomes
\begin{equation}\label{eqn.aug-PG}
	\mbox{select } \vek s_{j}\in\CU\mand\vek t_{j}\in\CV_{j}\msuchthat \vek b - \vek A\prn{\vek x_{0} + \vek s_{j}+\vek t_{j}}\perp \prn{\widetilde{\CU}+\widetilde{\CV}_{j}}.
\end{equation}
We represent these subspaces with the four matrices
$\vek U,\,\widetilde{\vek U} \in \Cmn{n}{k}$ and $ \vek V_{j},\,\widetilde{\vek V}_{j} \in \Cmn{n}{j}$, such that $\Ra{\MU} = \CU$, $\Ra{\widetilde{\vek U}}=\wtl{\CU}$, etc.
Associated to every such augmented projection method are a pair of projectors: $\vek P = \vek U\prn{\widetilde{\vek U}^{*}\vek A\vek U}^{-1}\widetilde{\vek U}^{*}\vek A$, the projector onto $\CU$ along $\prn{\vek A^{*}\widetilde{\vek U}}^{\perp}$, and $\vek Q = \vek A\vek U\prn{\widetilde{\vek U}^{*}\vek A\vek U}^{-1}\widetilde{\vek U}^{*}$, the projector onto $\vek A\,\CU$ along $\widetilde{\CU}^{\perp}$.
Analogous to the derivation of (\ref{eqn.derive-proj-Axb}), see also \cite{dSKS.2018}, the constraint
\cref{eqn.aug-PG} can be reformulated as the projected approximation problem whereby we take $\MV_j \vek y_{j}$ as the approximate solution of
\begin{equation}\label{eqn.projAxb}
	\prn{\vek I -\vek Q}\vek A \Vt = \prn{\vek I -\vek Q}\vek r_{0}
\end{equation}
that satisfies, for the residual of (\ref{eqn.projAxb}), $\wh{\Vr}_j$, the orthogonality condition
\begin{equation}\nn
	\wh{\vek r}_{j} = \prn{\vek I -\vek Q}\prn{\vek r_{0} - \vek A\vek V_{j}\vek y_{j}} \perp \widetilde{\CV}_{j} \quad \mbox{  and so  }\quad \vek t_{j}=\vek V_{j}\vek y_{j},
\end{equation}
and we use $\vek P$ to compute
\begin{equation}\nn
	\vek s_{j} =\vek P\Vt -\vek P\vek t_{j}
    \mwhere \vek P\vek \Vt = \vek U\prn{\widetilde{\vek U}^{\ast}\vek A\vek U}^{-1}\widetilde{\vek U}^{\ast}\vek r_{0},
\end{equation}
which (when we substitute in the expression for $\vek P$ and $\vek t_{j}=\vek V_{j}\vek y_{j}$)
leads to the construction of the full approximation at iteration $j$
\begin{equation}\label{eqn.general-approx-construction}
	\vek x_{j} = \vek U\prn{\widetilde{\vek U}^{\ast}\vek A\vek U}^{-1}\widetilde{\vek U}^{\ast}\vek r_{0} + \vek V_{j}\vek y_{j} - \vek U\vek B_{j}\vek y_{j}\qquad\mbox{where}\qquad\vek B_{j}=\prn{\widetilde{\vek U}^{*}\vek A\vek U}^{-1}\widetilde{\vek U}^{*}\vek A\vek V_{j}.
\end{equation}
We note that the matrix $\vek B_{j}$ is generally built column-by-column iteratively, and
that the residual of the full problem and the projected subproblem are, in fact, equal (so the residual norm is known without making the full update); i.e.,
\begin{equation}\label{eqn.proj-full-resid-eq}
	\vek r_{j} = \vek b - \vek A\prn{\vek x_{0} + \vek s_{j}+\vek t_{j}} = \widehat{\vek r}_{j}.
\end{equation}
This means properties of the projected coefficient matrix in \eqref{eqn.projAxb} and the iterative method we use to approximate its solution (i.e., our choices of
$\CV_{j}$ and $\widetilde{\CV}_{j}$) will dictate the convergence behavior of the augmented method.  Properties such as \eqref{eqn.general-approx-construction}
and \eqref{eqn.proj-full-resid-eq} are common to all subspace recycling methods.  These and other characteristics which can be gleaned from the framework
allow us to much more systematically design and implement a subspace recycling method.

\brem
Although one can choose any pair of subspaces
$\CV_{j},\,\widetilde{\CV}_{j}$, it is important that the choice results in an efficient method.
The review \cite{dSKS.2018} details how to leverage this theory to obtain customized recycling methods.  The greatest strength of viewing recycling methods through the lens of this framework is
that it decouples the choice of the iteratively generated correction and constraint spaces $\CV_{j}$ and $\widetilde{\CV}_{j}$ from the projected operator induced by
enforcing the residual constraint via \eqref{eqn.projAxb}.
Previously, these methods have been considered in a context where the Krylov subspace is generated by
the projected operator and an appropriate right-hand side.  However, this excludes useful augmentation schemes that can be described in this framework\cite{dSKS.2018}.
It is this observation which enables greater latitude in the systematic design of new customized recycling methods.

\erem

\paragraph{Example: Recycled FOM}
With the above framework and the techniques for building new methods thereof, it becomes much more
straightforward to build new, customized recycling methods using specific subspaces and residual constraints.  Following \cite{dSKS.2018}, we demonstrate this for a recycling FOM method.  We show two choices of projectors $\vek Q$ such that
$\CV_{j} = \CK_{j}\prn{\prn{\vek I - \vek Q}\vek A,\prn{\vek I - \vek Q}\vek r_{0}}$ leads to a viable implementation of a recycling FOM
method: in one case $\vek Q$ is the oblique projector onto
$\vek A\,\vek U$ along $\CU^{\perp}$, in the other case $\vek Q$ is the orthogonal projector onto $\vek A\,\CU$.

The first choice of $\vek Q$ is explored in detail to produce a full implementation.  In this case, we have
$\vek Q=\vek A\vek U\prn{\vek U^{\ast}\vek A\vek U}^{-1}\vek U^{\ast}$, and
as a consequence of \eqref{eqn.proj-relation},
$\vek P=\vek U\prn{\vek U^{\ast}\vek A\vek U}^{-1}\vek U^{\ast}\vek A$.
According to the framework, we  apply FOM to the projected
problem \eqref{eqn.proj-Axb} with the modification that,
during the Arnoldi process for generating a basis for
$\CK_{j}\prn{\prn{\vek I - \vek Q}\vek A,\prn{\vek I - \vek Q}\vek r_{0}}$,
we store the coefficients
$\vek B_{j}=\prn{\vek U^{\ast}\vek A\vek U}^{-1}\vek U^{\ast}\vek A\vek V_{j}$
which come from applying $\vek Q$ to $\vek A\vek v_{i}$ for each $i$.  Then
we solve the FOM linear problem
\begin{equation}\nonumber
	\vek H_{j}\vek y_{j} = \beta\vek e_{1},\qquad\beta=\nm{\prn{\vek I - \vek Q}\vek r_{0}},
\end{equation}
and at the end of the iteration we set $\vek x_{j} = (\vek x_{0} + \vek V_{j}\vek y_{j}) + \vek U\prn{\vek U^{\ast}\vek A\vek U}^{-1}\vek U^{\ast}\vek r_{0} - \vek U\vek B_{j}\vek y_{j}$.

\section{Practical realizations of the recycling framework}\label{sec.practical-realizations}
Although the framework described in \Cref{section.framework} can be used to understand the vast majority of subspace recycling methods,
most methods  were not derived this way.
There are some exceptions, notably those arising from
earlier proposed recycling frameworks \cite{Gaul.2014-phd,GGL.2013,Gutknecht.biCG-aug.2014} such as \cite{Gaul.S.rMINRES.2015} and the methods
proposed in the forthcoming review \cite{dSKS.2018}.
In this section, we
give an overview of existing methods and put them into context of the framework in \Cref{section.framework}.
We note that when discussing practical implementations of these methods, one must consider the developments in two research communities with
overlapping interests, those for
solving linear systems arising from discretizations of \emph{well-posed problems} (which may be ill-conditioned due
to isolated clusters of relatively small eigen- or singular values)
and those arising from discretizations of \emph{ill-posed problems} (which are characterized by rapidly decreasing singular values with no large, distinct gaps).
There are augmented/recycling-type methods arising from both communities, where similar methods often arose by happenstance almost in parallel.
\subsection{Full basis storage methods}
Methods for non-Hermitian systems, which do not take advantage of a short-recurrences, require that the basis (Arnoldi) vectors from each iteration
be stored until the end of the current restart cycle.  Thus, recycling for these methods must be between restart cycles for a system associated to
coefficient matrix $\vek A^{(i)}$ as well as between solves for $\vek A^{(i)}$ and $\vek A^{(i+1)}$.  Indeed, methods such as GMRES-DR
\cite{Morgan.GMRESDR.2002} can be understood as recycling exclusively between restart cycles using harmonic Ritz vectors.  Thus, as described in \Cref{sec.history-recyc}, many of these methods (see \Cref{sec.res-min-methods}) can be understood as the confluence of the
ideas of reuse of information between cycles and a more general passing of information between different linear systems in a systematic manner.

\subsubsection{Residual minimization methods}\label{sec.res-min-methods}
Although they are not often explicitly derived as such, recycling (or, more generally, augmentation) methods which minimize the residual
over the augmented space $\CU + \CV_{j}$ can be characterized with the well-known minimum residual constraint
\begin{equation}\nonumber
	\mbox{select } \vek s_{j}\in\CU\mand\vek t_{j}\in\CV_{j}\msuchthat \vek b - \vek A\prn{\vek x_{0} + \vek s_{j}+\vek t_{j}}\perp \vek A\,\prn{\CU+\CV_{j}}.
\end{equation}
which is a specific case of the more general (\ref{eqn.aug-PG}).  This is the theoretical umbrella which covers all such methods.  Two decisions then
determine exactly which method (up to implementation) is being proposed: which Krylov subspace is represented by $\CV_{j}$ and what method of subspace
downselection is being considered.
\paragraph{GCRO-DR/Recycled GMRES}
A strength of GCRO-based approaches was that they allowed one to minimize the residual over the combination of an arbitrary
subspace $\CU$ and a Krylov subspace generated from the projected coefficient matrix and right-hand side from (\ref{eqn.projAxb}).
In retrospect, this is one realization of the
notion of enforcing a constraint on an augmented subspace, as in \Cref{section.framework}. For all GCRO-based methods, we have the common
choice that $\CV_{j} = \CK_{j}\prn{\prn{\vek I-\vek Q} \vek A,\prn{\vek I-\vek Q}\vek r_{0}}$ with $\vek Q$ being the orthogonal projector
onto $\vek A\,\CU$.

Practically speaking, one generates via the Arnoldi process an orthonormal basis for the subspace
$\CK_{j}\prn{\prn{\vek I-\vek Q} \vek A,\prn{\vek I-\vek Q}\vek r_{0}}$. Let $\vek C\in\Cmn{n}{k}$, $\MC^*\MC=\MI$, and $\Ra{\MC} = \vek A\,\CU$ so
that $\vek Q=\vek C\vek C^{\ast}$.  With this, the Arnoldi relation becomes
\begin{equation}\label{eqn.proj-arnoldi}
	\prn{\vek I-\vek Q} \vek A\vek V_{j} = \vek V_{j+1}\underline{\vek H}_{j} \iff \vek A\vek V_{j} = \vek C\vek B_{j} + \vek V_{j+1}\underline{\vek H}_{j}\qquad\mbox{where}\qquad \vek B_{j}=\vek C^{\ast}\vek A\vek V_{j},
\end{equation}
and this yields the modified Arnoldi relation
\begin{equation}
	\vek A \begin{bmatrix} \vek U & \vek V_{j} \end{bmatrix} =  \begin{bmatrix}\vek C & \vek V_{j+1} \end{bmatrix} \begin{bmatrix} \vek I & \vek B_{j}\\  0 & \underline{\vek H}_{j}\end{bmatrix}.
\end{equation}
Using, in addition, $\Vr_0 = \MQ\Vr_0 + (\MI - \MQ)\Vr_0$ and
$(\MI - \MQ)\Vr_0 = \Vv_1 \nm{(\MI - \MQ)\Vr_0 }$, one can directly derive a practical implementation of the GCRO minimization yielding $\vek x_{j} = \vek x_{0} + \vek s_{j} + \vek t_{j}$ where,
\begin{equation}\label{eqn.gcro-ls-min}
	\prn{\vek z_{j}, \vek y_{j}} = \argmin{\vek z \in \Cn{k}\atop
\vek y \in\Cn{j}} \nm{\begin{bmatrix} \vek I & \vek B_{j}\\ 0 & \underline{\vek H}_{j}\end{bmatrix}
\begin{bmatrix} \vek z \\ \vek y \end{bmatrix} -
\begin{bmatrix}
  \MC^*\Vr_0 \\
  \Ve_1\beta
\end{bmatrix}} \mand \vek s_{j} = \vek U\vek z_{j} ,\,\vek t_{j} = \vek V_{j}\vek y_{j},\qquad \mbox{where } \beta = \nm{\prn{\vek I - \vek Q}  \vek r_{0}},
\end{equation}
where $\vek e_{1}\in\bbC^{k+j+1}$ is the $1$st Euclidean basis vector.
In \cite{EDS.GCRO.1996}, it was also suggested to  solve
\eqref{eqn.gcro-ls-min} blockwise, first for $\vek y_{j}$ and then
set $\vek z_{j} = \MC^*\Vr_0 - \MB_j \Vy_j$.
This is equivalent to applying GMRES directly to the projected problem (\ref{eqn.proj-Axb}) and storing in $\vek B_{j}$ the coefficients obtained from applying $\vek Q$ during the Arnoldi iteration (\ref{eqn.proj-arnoldi}).
This also falls directly out of the framework in \Cref{section.framework}:
one can work out that for initial error $\vek t$,
\begin{equation}\label{eqn.gcrodr-update}
	\vek x_{j} = \vek x_{0} + \vek P\vek t + \prn{\vek I-\vek P}\vek t_{j} \iff \vek x_{0} + \vek U\vek C^{\ast}\vek r_{0} + \vek V_{j}\vek y_{j} - \vek U\prn{\vek B_{j}\vek y_{j}}
\end{equation}
where $\vek y_{j}$ satisfies (\ref{eqn.gmres-ls-min}) for $\underline{\vek H}_{j}$ generated from (\ref{eqn.proj-arnoldi}) and
$\beta = \nm{\prn{\vek I - \vek Q}  \vek r_{0}}$.  The form of the update (\ref{eqn.gcrodr-update}) shows that this method can be implemented
as a GMRES iteration applied to an equation of the form (\ref{eqn.projAxb}) where we must simply store $\vek B_{j}$ as it is computed.  This enables
the construction of the full GCRO-DR approximation using the coefficients $\vek y_{j}$ obtained from applying GMRES to (\ref{eqn.projAxb}).
This is a specific realization of the updating formula (\ref{eqn.general-approx-construction}),
reinforcing the notion that all
residual constraint, augmentation-based subspace recycling methods share
this structure, which allows
them to be understood as known iterative methods being applied
to projected linear systems. Multiple authors have made this conclusion for specific subspace recycling algorithms \cite{KR.2012,GGL.2013,Gutknecht.biCG-aug.2014,Gaul.2014-phd,EDS.GCRO.1996}.

This setup is valid for any GCRO-based recycling method.
What differentiates most GCRO-based recycling methods is how we choose a subspace
of $\CU + \CV_{j}$ to
generate and $\CU_{new}$, the new recycled subspace.
The GCRO-T method \cite{deSturler.GCROT.1999} is focused on minimizing the penalty
for discarding some vectors at the end of a restart cycle
in a process called \emph{optimal truncation}.
For further details see \Cref{section.whatrecyc}.
A number of flexible and block variants have also been proposed in the
 literature \cite{HZ.simp-FGCROT.2010,MZL.2014,PSS.2016}.

The GCRO-DR method \cite{Parks.GCRODR.2006} builds on the same
 GCRO-type minimization but combined with the deflated-restarting strategies of Morgan \cite{Morgan.GMRESDR.2002},
wherein the subspace retained between restart cycles is taken to be some harmonic Ritz vectors.
There are strong associations between this strategy for solving linear systems and the implicitly restarted Arnoldi strategy for
the computation of eigenvalues and eigenvectors of large, sparse matrices.
Indeed, as was discussed in \Cref{sec.history-recyc}, if one looks at the forbearers of the *-DR
strategies, one sees that many elements take inspiration from the implicit restarting techniques for eigen-computations \cite{Morgan2000}.  We discuss this strategy further
in the context of its effectiveness as a recycling technique along in a wider discussion of recycling strategies in \Cref{section.whatrecyc}.
A block version of GMRES-DR has been proposed
\cite{Morg_2005} and the same is true for GCRO-DR \cite{PSS.2016}.

\paragraph{Recycled GMRES with an alternative projector}
As noted in\cite{GGL.2013,Gaul.2014-phd}, there is a second variant of augmented GMRES that fits into the augmentation
framework that was first proposed in \cite{Erl.Nabben.Defl-balPrec.2008} in the context of deflation-type GMRES methods.
This variant fits into the framework in \Cref{section.framework}.
The setup is similar to a GCRO-based method, except that the residual minimization is performed over a
obliquely projected Krylov subspace. 
This second variant was used to propose an alternative type of recycled MINRES
\cite{Gaul.S.rMINRES.2015} than that which is described in
\Cref{section.short-recur}.

In this case, let $\vek P$ be the oblique projector onto the subspace $\CU$ along $\prn{\vek A^{\ast}\,\CU}^{\perp}$ and $\vek Q$ be the sibling projector
satisfying (\ref{eqn.proj-relation}), the oblique projector onto $\vek A \CU$ along $\CU^{\perp}$.  We use $\vek P$ to split the initial error
$\vek t = \vek P\vek t + \prn{\vek I - \vek P}\vek t$, observing that as before $\vek P\vek t$ can be computed explicitly.
We then approximate $\prn{\vek I - \vek P}\vek t$ by applying GMRES to (\ref{eqn.proj-Axb}) obtaining $\vek t_{j}$.  The full approximation is then
$\vek x_{j} = \vek x_{0} + \vek P\vek t + \prn{\vek I-\vek P}\vek t_{j}$.  This is equivalent to minimizing the full residual over
the subspace $\CU + \CV_{j}$ where \linebreak $\CV_{j} = \CK_{j}\prn{\prn{\vek I - \vek Q}\vek A,\prn{\vek I - \vek Q}\vek r_{0}}$ where $\vek Q$ is an
oblique projector rather than an orthogonal projector, as in GCRO.

\paragraph{Augmented GMRES for ill-posed problems}
In \cite{BR.2007}, the authors propose to reconstruct the solution of a non-Hermitian ill-posed problem of the form
\eqref{eqn.Axb} over an augmented Krylov subspace as well as over a preselected space.
The authors propose to do this by decomposing
the problem into two subproblems using projectors.  The larger subproblem is posed in a space orthogonal to the augmenting
subspace and is solved using an iterative method, and the
smaller problem is posed in the augmenting subspace and is solved directly.
This strategy can be cast as an augmented method\cite{BR-2.2007}.
Indeed, this can be directly related to the technique of splitting the error using a projector (\ref{eqn.error-splitting}). 
Applying GMRES to the projected subproblem leads to a method which is mathematically equivalent to GMRES with recycling,
i.e., setting $\CV_{j} = \CK_{j}\prn{\prn{\vek I -\vek Q}\vek A,\prn{\vek I -\vek Q}\vek r_{0}}$.
However, in this context, other choices of Krylov subspaces were suggested.  Often, it has been shown beneficial
to employ a strategy known as range-restriction, wherein one
uses a power of the operator times the right-hand side, e.g.,
$\CV_{j} = \CK_{j}\prn{\prn{\vek I -\vek Q}\vek A,\prn{\vek I -\vek Q}\vek A\vek r_{0}}$.
We note that based on the work in \cite{BR-2.2007},
another adaptive augmented method has also been developed \cite{KN.aug-GMRES-ip}.
This body of literature is not necessarily concerned with recycling spectral information as much as it is with an augmentation space that encodes certain known features of the solution.
However, strategies advocated in \cite{ChuStuJia_20} show that
recycling in the sense of GCRO-DR \cite{Parks.GCRODR.2006} can also be effective.

\paragraph{Unprojected GMRES with range-restriction}
The fact that the modified Arnoldi iteration from \cite{BR-2.2007} indeed is equivalent to a Krylov subspace iteration
for the composition of the operator $\vek A$ and a projector was observed in \cite[Remark 2.1]{DGH.2014}, and this
remark leads the authors to propose a modification of the augmentation method\cite{BR-2.2007}.  The authors assert that augmenting an unprojected Krylov
subspace allows a residual polynomial approximating the polynomial
representation of the true inverse to be constructed.  	
Using the polynomial representation of the inverse may not be the best approach
for considering the effectiveness of this method, but we can evaluate this augmentation strategy nonetheless for various situations. They also assert
that a poor choice of $\CU$ causes the iteration to go awry with no
chance of recovery.
They propose reordering the steps of the modified
Arnoldi process such that a GMRES-type iteration does reconstruct part of the approximation over a Krylov subspace
generated by $\vek A$ rather than the projected operator of the form in (\ref{eqn.proj-Axb}).
We note that with the reordering of orthogonalization steps, this method is actually closely related to the flexible GMRES-based augmentation
scheme of \cite{Chapman1997} discussed in \Cref{sec.history-schemes}.

Rather than generating the residual
$\prn{\vek I - \vek Q}\vek r_{0}$ and generating a Krylov subspace with respect to the projected operator,
in \cite{DGH.2014}, the
authors propose to generate $\CK_{j}\prn{\vek A,\vek p_{0}}$, where
$\vek p_{0}\in\curl{\vek r_{0},\vek A\vek r_{0}}$
\footnote{The rationale for choosing $\vek p_{0}=\vek A\vek r_{0}$ is that the right-hand
side or initial residual for an ill-posed problem is noisy and can corrupt the iterative solution process.  Since $\vek A$ is assumed to be the discretization
of an operator which has smoothing properties,
$\vek A\vek r_{0}$ is smoother and will contain less noise.  In Hanke's monograph \cite{H1995-book} on regularization properties
of various iterative methods, this was denoted in the context of MINRES by MR2.}.
After each new Arnoldi vector is generated, it is used to project the image of the augmentation
space under the action of the operator away from that vector.  In other words, at iteration $j$, the algorithm generates
$\vek W_{j} = \prn{\vek I - \vek V_{j+1}\vek V_{j+1}^{\ast}}\vek A\vek U$ progressively.
This leads to
\begin{equation}\nonumber
	\vek A\begin{bmatrix} \vek V_{j} & \vek U \end{bmatrix} = \begin{bmatrix} \vek V_{j+1} & \vek W_{j}  \end{bmatrix}\begin{bmatrix} \underline{\vek H}_{j} & \vek L_{j}\\ 0 & \vek F_{j} \end{bmatrix},\mwith \vek L_{j} = \vek V_{j+1}^{\ast}\vek C\mand \vek F_{j} =\vek W_{j} ^{\ast}\vek C.
\end{equation}
This leads to a GMRES-like iteration in which one must solve a small least squares problem
of the form similar to \eqref{eqn.gmres-ls-min}
\begin{equation}\nonumber
	\begin{bmatrix}\vek y_{j}\\ \vek z_{j}\end{bmatrix} = \argmin{\vek z\in\Cn{k}\atop \vek y\in\Cn{j}}\nm{\begin{bmatrix} \underline{\vek H}_{j} & \vek L_{j}\\ 0 & \vek F_{j} \end{bmatrix}\begin{bmatrix}\vek y\\ \vek z\end{bmatrix} - \begin{bmatrix} \vek V_{j+1} & \vek W_{j}  \end{bmatrix}^{\ast}\vek r_{0}}, \qquad\mbox{where}\qquad \vek s_{j} = \vek U\vek z_{j}\mbox{  and  }\vek t_{j} = \vek V_{j}\vek y_{j}.
\end{equation}
They call this method Regularized range-restricted GMRES (R3GMRES).
It should be noted that this strategy also fits into the general augmentation framework, with $\CV_{j}= \CK_{j}\prn{\vek A,\vek p_{0}}$.
Further details and elaboration on the relationship of this method to the framework in
\Cref{section.framework} can be found in the upcoming paper \cite{KMS-aug-arn-tikh.2020}.

\subsection{Short-recurrence-based methods}\label{section.short-recur}
If the coefficient matrices in (\ref{eqn.Aixibi}) are Hermitian or one uses a short recurrence method for non-Hermitian systems,
the additional challenges that arise for subspace recycling are different from those for
GMRES and other full basis storage methods. Short recurrence methods do not need to restart, so there is no need to recycle for a subsequent restart cycle.
The additional challenge is determining how to
downselect the constructed basis to recycle.  There is no need to recycle at a restart for the current system, but one should select a recycle space for the subsequent problem.  The usual strategy is to store a running window of the $p$ most recent (Hermitian or bi-orthogonal) Lanczos vectors as columns of a matrix $\MV_{curr}\in\Cmn{n}{p}$.
Consider that we have a recycled subspace $\CU_{curr}$ being used for the current system.  We then have a separate recycled subspace $\CU_{next}$ which is held for the next system.  We initialize $\CU_{next}\leftarrow\CU_{curr}$.
When $p$ vectors have been
stored, the existing recycled subspace $\CU_{next}$ is overwritten by computing a subspace of $\CU_{next} + {\rm range}\prn{\MV_{curr}}$ according to the chosen downselection criteria.   The vectors in $\MV_{curr}$ are discarded and
the matrix is filled again with the next set of $p$ Lanczos vectors.  Examples are discussed in \cite{Ahuja.Recyc-BiCG.2012,KilmerdeSturler2006,Wang.TopOptRecyc.2007,MR1270124,Morgan.rest-Lanczos-defl.2010};
a variant of this idea is discussed in \cite{MR2599785}.

\paragraph{$\vek A$-norm optimal error methods for Hermitian positive definite systems}
Augmented conjugate gradient-type methods have been proposed independently in both the well- and ill-posed problems communities.  Discretizations of
many elliptic and parabolic partial differential equations lead to HPD discrete linear problems.  For discrete ill-posed problems, one often uses
a regularization of  the normal equations
associated to the linear problem, which generally produces a Hermitian positive-definite problem.
\footnote{It has been shown that CG applied
to this problem with an appropriate early stopping rule satisfies the formal definition of regularization; see \cite{EHN.1996-book} and \cite{H1995-book}
for details.  Indeed, augmentation-type methods fitting into the framework in \Cref{section.framework} have been shown to also satisfy the formal definition
of a regularization in \cite{RRKMS.regrecyc.2020}.}
An augmented conjugate gradients approach was first proposed in 2000 in \cite{SYEG.2000}.
This
was then followed up by \cite{ErhGuy_2000} wherein the method was improved and it was used to treat a sequence of systems with the same HPD linear system but changing right-hand sides.  Some of the residuals from previous systems are proposed to span the augmentation space $\CU$.
An unpublished manuscript \cite{PSS.RCG.2016} proposes a recycled variant of the augmented CG method that employs the strategy above, storing a fixed window of $\vek A$-conjugate directions to update the recycled subspace for the next system.
In the discrete ill-posed problem setting, a similar method was proposed in \cite{CRS.2003}.  However, the authors go further by modifying the minimization
over the augmented space to transform it into a Tikhonov-type penalized minimization.  An implicitly restarted Lanczos method for symmetric eigenvalue
problems also has elements fitting into this framework \cite{MR1270124}.

The CG algorithm can be formulated as a highly efficient algorithm by exploiting the fact that an HPD operator can define an inner product. It is important to exploit this
property of the linear system also in the recycling version.
Several efficient recycling-like algorithms have been proposed
doing this \cite{SYEG.2000,CRS.2003,KR.2012,Gaul.2014-phd}.
While the framework above allows us to define different recycling methods by
distinct choices for the recycle correction and constraint spaces and inner
product, we obtain a very efficient method by exploiting
the $\MA$-inner product and using the Galerkin approach, that is, the
constraint space equals the correction space. This leads to
$\MQ = \MA \MU (\MU^*\MA\MU)^{-1}\MU^*$ and,
following (\ref{eqn.proj-relation}),
$\MP = \MU (\MU^*\MA\MU)^{-1}\MU^*\MA $, which
leads to several special properties.
By inspection, $\MQ = \MP^*$ and
$(\MI - \MQ)\MA = \MA(\MI - \MP)$ is Hermitian
(the equality holds for any $\MQ$ and $\MP$
satisfying (\ref{eqn.proj-relation})).
Finally, we have that $(\MI - \MQ)\MA = \MA(\MI - \MP)$ is
positive semi-definite, as $\MA$ defines an inner product:
$\Vw^*(\MI - \MQ)\MA\Vw = \Vw^*(\MI - \MQ)(\MI - \MQ)\MA\Vw =
\Vw^*(\MI - \MP)^*\MA(\MI - \MP)\Vw \geq 0$ and
$\Vw^*(\MI - \MQ)\MA\Vw = 0 \eqv (\MI - \MP)\Vw = \Vo$.

Without loss of generality, we consider henceforth the well-posed problem with $\vek A$ HPD and $\MQ$ and $\MP$ as defined above.
Since (\ref{eqn.proj-Axb}) is consistent, we can apply CG directly to this problem.
In fact, since $\Ra{(\MI-\MQ)} = \CU^{\perp}$,  $\CK_j((\MI-\MQ)\MA,(\MI-\MQ)\Vr_0) \subseteq \CU^{\perp}$ and
hence for any $\Vw \in \CK_j((\MI-\MQ)\MA,(\MI-\MQ)\Vr_0)$, $(\MI - \MP)\Vw \neq \V0$.
Taking $\Vv_1 = \b^{-1} (\MI - \MQ)\Vr_0$ with
$\b = \| (\MI - \MQ)\Vr_0 \|$ the Lanczos relation becomes
(cf. \cref{eqn.Lanc_rec-2})
\[
  (\MI - \MQ)\vek A \vek V_j  =  \vek V_{j+1}\underline{\vek T}_j =
  \vek V_j \vek T_j + \vek v_{j+1} h_{j+1,j} \Ve_{j}^*.
\]
Following the above, $\MT_j$ is positive definite and therefore
the LU decomposition,
$\MT_j = \ML_j \MU_j$ (without) pivoting exists. So, we can again
apply the change of basis transformation $\MW_j = \MV_j  \MU_j^{-1}$,
set $\wtl{\Vy}_j = \ML_j^{-1} \Ve_1$, and run the standard CG iteration on
the projected system. The (full) solution can then be computed according to
\eqref{eqn.general-approx-construction}. In this case, by construction
$\Vr_j \perp  \{\Vu_1, \ldots, \Vu_k, \Vv_1, \ldots, \Vv_j \}$, which proves
that the error is minimized in the $\MA$-norm over the space $\Ra{\MU} \oplus \Ra{\MV}_j$.

In some alternative approaches \cite{SYEG.2000,CRS.2003,KR.2012,Gaul.2014-phd}, a change of basis is used that generates an $\MA$-orthogonal basis for
$\Ra{\MU} \oplus \Ra{\MV}_j$, which could be advantageous for some
applications. Other alternative approaches compute more
accurate eigenvectors over multiple linear systems, but
for each linear system, they deflate these only from the initial
residual \cite{MR2599785,Morgan.rest-Lanczos-defl.2010}.

\paragraph{Minimum residual methods for Hermitian indefinite systems}
Consider (\ref{eqn.Aixibi}) with $\MA$ Hermitian and indefinite, initial guess
$\Vx_0$, and residual $\Vr_0$. As in \Cref{section.framework}, let
$\MU = [\Vu_1, \, \Vu_2, \, \ldots,\, \Vu_k]$ define the recycle correction space
$\CU = \Ra{\MU}$, and $\MA \MU = \MC$ with
$\MC^* \MC = \MI$. The recycle constraint space is chosen to be
$\wtl{\CU} = \Ra{\MC}$. Hence, $\MQ = \MC \MC^*$
and $\MP = \MU \MC^* \MA$, which gives $\MP\Vt = \MU \MC^* \Vr_0$.
Recycling MINRES (rMINRES)\cite{Wang.TopOptRecyc.2007} works by applying MINRES to approximately solve the system
$(\MI - \MQ)\MA \Vt = (\MI - \MQ)\Vr_0$ for an approximation to the update $\Vt$. Since
$\CK((\MI - \MQ)\MA,(\MI - \MQ)\Vr_0) \subset \Ra{\MI - \MQ}$,
$(\MI-\MQ)\MA \! : \! \Ra{\MI - \MQ} \rite \Ra{\MI - \MQ} $, and $\MA = \MA^*$,  we have for
all $\Vz_1, \Vz_2 \in \Ra{\MI - \MQ}$,
\[
  \la (\MI - \MQ)\MA \Vz_1, \Vz_2\ra =
  \la (\MI - \MQ)\MA(\MI - \MQ) \Vz_1, \Vz_2\ra =
  \la \Vz_1, (\MI - \MQ)\MA(\MI - \MQ)\Vz_2\ra =
  \la \Vz_1, (\MI - \MQ)\MA\Vz_2\ra .
\]
So, $(\MI - \MQ)\MA$ is self-adjoint (`Hermitian') over
$\CK((\MI - \MQ)\MA,(\MI - \MQ)\Vr_0) $, and we use the
Hermitian Lanczos process (\ref{eqn.Lanc_rec-1}) to get
\eqs \label{eqn.proj-lanczos_1}
  (\MI - \MQ)\MA \MV_j & = & \MV_{j+1}\ul{\MT}_j \quad \eqv \quad
      \MA\MV_j = \MC \MB_j + \MV_{j+1}\ul{\MT}_j ,
\eqe
with $\MB_j = \MC^*\MA\MV_j$.
As for standard MINRES, we can apply a change of basis,
$\MW_j \MRr_j= \MV_j$, using the thin QR decomposition
$\ul{\MT}_j = \ul{\MG}_j \MRr_j$,
set $\wtl{\Vy} = \ul{\MG}_j^*\Ve_1\|(\MI - \MQ)\Vr_0\|$,
and during the iteration use the (partial) solution update
$\Vx_j = \Vx_{j-1} + \Vw_j \wtl{y}_j  =  \Vx_0 + \MW_j \wtl{\Vy}_j$,
where $\Vt_j = \MW_j \wtl{\Vy}_j$ and $\wtl{y}_j$ is the $j$-th component
of $\wtl{\Vy}_j$.
At the end of the iteration, we set
\eqs \label{eqn.rminres-update}
  \Vx_j & = & \Vx_j + \MP \Vt - \MP \Vt_j =
    \Vx_0 + \MW_j \wtl{\Vy}_j + \MU\MC^*\Vr_0 -  \MU \MB_j (\MRr^{-1} \wtl{\Vy}_j ).
\eqe
This approach postpones all $\MU$ updates during the iteration to a
single update at the end, which saves $O(kn)$ work per iteration \cite{Wang.TopOptRecyc.2007}.
Several other efficiency improvements
are discussed in \cite{Wang.TopOptRecyc.2007,Recyc-Imp-Tom.2010}. This
includes a change of basis that allows to discard the columns of $\MB_j$ as
we go (possibly important if many iterations are required) and very efficient recurrences to compute a recycle space.

\paragraph{Bi-orthogonal Lanczos-based methods and Transpose-free Variants}
For non-Hermitian system matrices, we consider BiCG and BiCGstab-based methods\cite{A.2009,Ahuja.Recyc-BiCG.2012,Ahuja.Sturler.rGCROT-rBiCGStab-Hybrid.2015,Ahuja.Recyc-BiCGStab-MOD.2015}.
Obviously, recycling versions of QMR \cite{Freund.QMR} and TFQMR \cite{Freund.TFQMR}
can developed as well. Recently, a recycling IDR(s) variant was developed \cite{NeuGreif_18}.

\paragraph{Recycling BiCG}
We modify the BiCG algorithm to use recycle spaces.
Here, we follow the approach chosen in \cite{Ahuja.Recyc-BiCG.2012}.
Let $\Ra{\MU} = \CU$ be the chosen recycle (correction) space for the primary linear system (\ref{eqn.Aixibi}) and $\Ra{\wh{\MU}} = \wh{\CU}$ be the chosen recycle correction space for the dual system.
Moreover,
we choose $\MU$ and $\wh{\MU}$ such that
$\MC = \MA \MU$ and $\wh{\MC} = \MA^* \wh{\MU}$
satisfy $\wh{\MC}^*\MC = \MD_c$ is real, diagonal, and invertible, i.e.,
$\wh{\MC}$ and $\MC$ are bi-orthogonal. This can always be done (for example
using the SVD \cite{Ahuja.Recyc-BiCG.2012} or an LDU decomposition), although
it may require reducing the dimension of the recycle spaces if
$\MD_c$ had one or more zeros on the diagonal.
We define $\MQ = \MC \MD_c^{-1} \wh{\MC}^*$ and
$\MP = \MU \MD_c^{-1} \wh{\MC}^*\MA$.
Following the discussion on the non-Hermitian Lanczos process in section \ref{section.kryl-basic}, we build dual bases for the subspaces
$\CK_{j} ((\MI-\MQ)\MA,(\MI-\MQ)\Vr_{0})$ and $\CK_{j}((\MI - \MQ^*)\MA^*,(\MI - \MQ^*)\wh{\Vr}_{0})$, where $\wh{\vek r}_{0}$ is either the initial residual of a dual problem involving $\vek A^{\ast}$, $\vek r_{0}$ itself, or some other non-zero vector.
This implies that the recycle constraint space is given by
$\wtl{\CU} = \Ra{\wh{\MC}}$, and the recycled constraint space for the
dual problem, if defined, is $\Ra{\MC}$.
Note that
$(\MI-\MQ)\MA: \Ra{\wh{\MC}}^{\perp} \rite \Ra{\wh{\MC}}^{\perp}$ and
$(\MI-\MQ^*)\MA^* : \Ra{\MC}^{\perp} \rite \Ra{\MC}^{\perp}$,
and we have
for any $\Vw \in \Ra{\wh{\MC}}^{\perp}, \wh{\Vw} \in \Ra{\MC}^{\perp}$
that
\eqs \nonumber
  \la (\MI - \MQ)\MA \Vw , \wh{\Vw} \ra & = &
  \la (\MI - \MQ)\MA(\MI - \MQ) \Vw , \wh{\Vw} \ra =
  \la \Vw ,   (\MI - \MQ^*)\MA^*(\MI - \MQ^*) \wh{\Vw} \ra
  =   \la \Vw ,  (\MI - \MQ^*)\MA^* \wh{\Vw} \ra.
\eqe
So, $(\MI - \MQ)\MA $ and $(\MI - \MQ^*)\MA^*$ act as adjoints
over the primary and dual Krylov spaces, and we can again apply a coupled three-term recurrence as in
\Cref{section.kryl-basic}:
\begin{eqnarray}\label{augBiLanOld}
	(\MI - \MQ)\MA \MV_j =
        \MV_{j+1}\underline{\MT}_j
    & \mbox{   and   } &
	(\MI - \MQ^*)\MA^* \wh{\MV}_j =
        \wh{\MV}_{j+1}\wh{\underline{\MT}}_j .
\end{eqnarray}

As for the standard BiCG discussion, we consider here the updates for the primary system; the updates for the dual system (if needed) can be computed analogously.
We compute (if it exists) the
LU-decomposition $\MT_j = \ML_j \MU_j$,
use the change of basis $\MW_j = \MV_j \MU_j^{-1}$, and set $\wtl{\Vy}_j = \ML_j^{-1} \Ve_1 \z$,
where $\z = \|(\MI - \MQ)\Vr_0\|$.
This allows us to eliminate the Lanczos vectors and update
\[
  \Vt_j = \Vt_{j-1} + \Vw_j \wtl{y}_j ,
\]
where $\wtl{y}_j$ is the $j$-th component of $\wtl{\Vy}_j$ (the new component).
The corresponding full solution is given by
\eqsn
  \Vx_j & = & \Vx_0 + \Vt_j + \MP \Vt - \MP \Vt_j  =
      \Vx_0 + \Vt_j + \MU \MD_c^{-1} \wh{\MC}^*\Vr_0 -
      \MU \MD_c^{-1} \wh{\MC}^*\MA \MW_j \wtl{\Vy}_j ,
\eqen
Only $\Vt_j$ is updated during the iteration.
The updates in the $\MU$ direction
are done after the final iteration, while the vector
$\MD_c^{-1} \wh{\MC}^*\MA\MW_j \wtl{\Vy}_j$ can be updated
during the iteration (without computing $\MA \MW_j$) \cite{Ahuja.Recyc-BiCG.2012}.
The matrices for the generalized eigenvalue problem that defines these subspaces can be constructed efficiently using recurrences\cite{Ahuja.Recyc-BiCG.2012}.

One can develop a recycling BiCGStab based on recycling BiCG using the polynomials defining the iteration vectors following \cite{vandervorst.bicgstab.1992};
see \cite{Ahuja.Recyc-BiCGStab-MOD.2015}. This leads to a recycling BiCGStab with the matrix $(\MI - \MQ)\MA$ and $\MQ$ as defined above. However, one can also take $\MQ$ as for recycling GCROT above, based on a single recycling
correction space $\CU$, yielding $\MQ = \MC\MC^*$; see \cite{Ahuja.Recyc-BiCGStab-MOD.2015} and \cite{Ahuja.Sturler.rGCROT-rBiCGStab-Hybrid.2015}, which includes an application where this approach is particularly useful.

\paragraph{LSQR-based methods}
This survey has focused on recycling methods for linear problems with square coefficient matrices.  However, the framework of
\Cref{section.framework} is not restricted to square matrices
or even matrix equations.  It is shown in \cite{dSKS.2018} that we can extend
this framework to Hilbert space operator equations $Tx=y$ where $T:\CX\rightarrow\CY$ is a linear mapping between two abstract Hilbert spaces.
One can approximate the solution using a Petrov-Galerkin residual constraint posed in $\CY$, with the correction spaces
$\CU,\,\CV_{j}\subset\CX$ and the constraint spaces $\widetilde{\CU},\,\widetilde{\CV}_{j}\subset\CY$.
This Hilbert space framework is also used in \cite{RRKMS.regrecyc.2020} to prove that augmentation methods
satisfy the formal definition of a regularization method in the infinite-dimensional ill-posed problems setting.

One realization of this more general notion of recycling methods arises when we have a linear system with a tall, skinny matrix wherein
we are seeking the least-squares solution.  Consider the (possibly inconsistent) linear system
\begin{equation}\nonumber
	\vek G\vek x \approx \vek f\ \  \mbox{where} \ \  \vek G\in\R^{m\times n}\ \  \mbox{and} \ \ m > n
\end{equation}
The LSQR method is a short-recurrence Golub-Kahan bidiagonalization (GKB) method \cite{PS-LSQR} which progressively
solves
\begin{equation}\nonumber
	\vek t_{j} = \argmin{\vek t\in\CK_{j}\prn{\vek G^{\ast}\vek G, \vek G^{\ast}\vek r_{0}} }\nm{\vek b - \vek G\prn{\vek x_{0} + \vek t}}\ \  \mbox{where} \ \ \vek r_{0} = \vek f - \vek G\vek x_{0},
\end{equation}
which can be shown to be equivalent to the residual constraint formulation
\begin{equation}\nonumber
		\mbox{select } \vek t_{j}\in\CK_{j}\prn{\vek G^{\ast}\vek G, \vek G^{\ast}\vek r_{0}}\msuchthat \vek f - \vek G\prn{\vek x_{0} + \vek t_{j}}\perp \vek G\vek G^{\ast}\CK_{j}\prn{\vek G\vek G^{\ast},\vek r_{0}}.
\end{equation}
The efficient progressive formulation of the algorithm arises from the fact that it is possible to generate simultaneously
via short recurrences orthonormal bases for the spaces $\CK_{j}\prn{\vek G^{\ast}\vek G, \vek G^{\ast}\vek r_{0}}$ and
$\CK_{j}\prn{\vek G\vek G^{\ast},\vek r_{0}}$ via the GKB.

With this formulation, one sees that it is possible to impose residual correction/constraint conditions over sums of spaces
for this problem just as in the square problem case (\ref{eqn.aug-PG}).  The main challenge then is to choose projected Krylov subspaces $\CV_{j}$
and $\widetilde{\CV}_{j}$ for the iteratively generated parts of the correction and constraint
spaces, respectively, which are related such that orthonormal bases can be generated via GKB-type short recurrences.
This has not been explored extensively in the literature, to our knowledge.  For acceleration of convergence for a single problem, a deflated-restart-type
method based on the theory presented for GMRES-DR \cite{Morgan.GMRESDR.2002} was proposed \cite{BRR-augLSQR.2013} .  A forthcoming paper for
sequences of regularized least-squares problems explores one particular strategy for choosing the projected Krylov spaces $\CV_{j}$ and $\widetilde{\CV}_{j}$
appropriately to ensure that the orthonormal bases can be generated via the GKB
\cite{ChuStuJia_20} and is further explored in the upcoming
review \cite{dSKS.2018}. 

\section{What to recycle}\label{section.whatrecyc}
It is difficult to make
general prescriptions about which subspaces to recycle because this depends on so many factors connected to the specific problem being solved.
We break down the different choices proposed
currently in the literature: approximate eigenvector augmentation, POD-type strategies, optimal truncation, and spaces from approximate solutions.

\paragraph{Approximate solution augmentation}
In \cite{KilmerdeSturler2006} solutions from previous nonlinear iterations are recycled to obtain good initial guesses in subsequent problems.
The success of this approach is problem dependent. If
the right hand sides do not change (much) from one linear system to the next
(as is the case in \cite{KilmerdeSturler2006}) and the
coefficient matrices are close, this typically produces good initial guesses.

\paragraph{Approximate eigenvector augmentation}
As discussed in \Cref{sec.history-recyc}, there has been limited analysis up to now on \emph{approximate} eigenvector augmentation
(see, e.g., \cite{Gaul.2014-phd,desturler.recycling-convergence-abstract.2012,KilmerdeSturler2006,morgan2020twogrid}
 and our forthcoming review \cite{dSKS.2018}), but there is strong evidence of the effectiveness of this strategy based on
actual results for application problems.  As these methods have often been shown to be equivalent to a deflation strategy, analysis pertaining to exact eigenvector deflation is also useful to study.

One rationale for augmenting with approximate eigenspaces is connected to the convergence theory for Krylov subspace methods.  If $\CU$ is an approximate
eigenspace, then the projector $\prn{\vek I - \vek Q}$ is an (in this case orthogonal) projector onto $\prn{\vek A\,\CU}^{\perp}$, and iterations during
a cycle of GCRO-DR take place in $\prn{\vek A\,\CU}^{\perp}$.  Thus, the projector $\prn{\vek I - \vek Q}$
may have the effect of damping possible negative influence on convergence speed of the approximated invariant subspace (following from the
theory introduced in \cite{Simoncini2005}), leading to
an accelerated convergence.  Care must be taken, however, as it has been shown that residual convergence need not necessarily be connected
to the spectral properties of the coefficient matrix \cite{Greenbaum.Any-Curve-Possible-GMRES.1996}.

Additionally, if two matrices $\vek A^{(i)}$ and $\vek A^{(i+1)}$  are ``close enough'',
then particular respective invariant subspaces may also be close.
For some differential operators (such as an elliptic operator),
the higher-frequency eigenmodes are associated with the larger eigenvalues, and this property usually is inherited by the discretized
matrix.
If the sequence of linear systems is induced by local changes in the matrix entries, then
the changes to the invariant subspaces associated to the higher-frequency eigenvectors (with larger eigenvalues) dominate.
The details are quite technical; so we direct the reader to \cite{KilmerdeSturler2006} where this was discussed
and quantified for the problem under consideration.
Conversely, for certain integral operators, one would want to recycle the approximate eigenvectors associated to
the larger Ritz values, again to capture the low-frequency eigenvectors.

\paragraph{Approximate singular vector augmentation}
In \cite{al2018recycling}, it has been recently proposed that one can also observe convergence speedups by recycling
 approximate singular vectors rather than eigenvectors, based on some analysis in \cite{Simoncini2000} that the residual having
 large components in the left singular vectors can cause reduced convergence speed.  The authors thus propose to use a Ritz-type approximation
of the left singular vectors (i.e., eigenvectors of $\vek A^{\ast}\vek A$) and recycle some of them.

\paragraph{Optimal truncation}
In \cite{deSturler.GCROT.1999}, a method is proposed to mitigate the effects of restarting after a cycle of GMRES, wherein the entire Krylov subspace
from the previous cycle is discarded.  Hence, we are disregarding orthogonality with respect to the discarded subspace.  This
causes the characteristic reduced convergence rate one sees with restarted GMRES as compared to full GMRES. In \cite{deSturler.GCROT.1999}, the author
develops a model to characterize this delay (called ``residual error'' in the paper).  One can again hearken back to the difference between the steepest descent method and
CG as in \Cref{sec.history-recyc} to understand this strategy.  Restarting necessitates ignoring orthogonality with respect to the previously generated search space, and one way to interpret the so-called residual
error caused by this is that the next cycle of GMRES induces a minimization over the new Krylov subspace which, since it is not orthogonalized against the previous, undoes some of the
improvement in the residual along directions from the previous Krylov subspace.  At worst, this can cause total stagnation and
it generally is known to cause a slowdown in convergence.  The GCROT method seeks to mitigate this problem with the described strategy, by maintaining orthogonality to some portion of the previously generated
Krylov subspace.

This is achieved via the assumption that subspaces which were important (for speed of convergence) to maintain orthogonality against will continue
to be important.  Thus, for a cycle with length $m$, one studies the convergence during the last $s<m$ iterations of the cycle and compares that to the
slower convergence which would have occurred had a dimension $k$ subspace been neglected during the Arnoldi orthogonalization in those last
$s$ iterations.
One can determine the dimension $k$ subspace that would have caused the most delay had orthogonality against that space been neglected\cite{deSturler.GCROT.1999}.  This subspace is then chosen to be retained for the next cycle in a process called \emph{optimal truncation}.
The model works just as well if we are selecting a subspace
of an augmented Krylov subspace of the form $\CU+\CV_{j}$, where
$\CV_{j}$ is a projected Krylov subspace of the form
$\CK_{j}\prn{\prn{\vek I-\vek Q}\vek A,\prn{\vek I-\vek Q}\vek r_{0}}$. To recycle between consecutive linear systems with respective coefficient matrices $\vek A^{(i)}$
and $\vek A^{(i+1)}$ one can also employ this strategy: if consecutive systems are ``close enough'', typically optimal truncation
will still confer benefits associated with maintaining orthogonality to the truncated subspace.

\paragraph{Proper Orthogonal Decomposition (POD)-type strategies}
In the setting of recycled GMRES for a non-Hermitian family of
shifted systems, the notion of using a POD-type strategy has actually been alluded to in one
of the numerical experiments \cite[Section 6.6]{S.2016}.  However, the author did not call this a POD strategy, and
this was not expanded upon.
This strategy takes advantage of the fact that for the test examples (coming from lattice quantum-chromodynamics
application problems) the solutions for all shifted systems with shifts in a positive interval suitably away from zero
depend smoothly on the shift.  Thus according to the theory of Kressner and Tolber \cite{KT.2011}, solutions corresponding to all
shifts in this positive interval can be well-approximated in a low-dimensional subspace; cf. \Cref{section.cont-param-dep} for more details.
Thus, a few solutions can be found, and the subspace they span can be used as an effective recycled subspace for solving the rest of the problems.

Specifically, a recycled CG method
was proposed in which a subspace of the correction space is selected for retention using a strategy based on techniques
from model order reduction \cite{Carlberg2015}. If one is solving a sequence of HPD problems,
and $j-1$ systems have already been solved, 
one should collect specially chosen snapshot vectors
accumulated from the first $j-1$ system solves and use them to generate the POD-subspace with which to augment
when applying CG to system $j$.  This subspace is generated by approximately minimizing the operator-norm-distance
between the true solution and its projection onto the POD subspace. 

\section{Exploiting operator structure and multiple right-hand sides}\label{section.additional-structure}
Sometimes, the sequence of coefficient matrices $\vek A^{(i)}$ have additional structure or for each $i$ we actually have a parameterized set of matrices $\vek A^{(i)}\prn{s}$,
each of which is associated to a linear system.  Furthermore, there may be many right-hand sides for a single system matrix, or slightly different right-hand sides for a sequence of related matrices.  Such problems present additional challenges for adapting
recycling.  Most often, this takes the form of families of matrices which have some linear shifting structure, although we do consider more general parameter dependence in \Cref{section.cont-param-dep}.
The most general formulation of a family of linear systems with linear shift structure is
\begin{equation} \label{eq:MRHS_all} (\MA^{(i)} + \gamma_\ell \ME) \Vx^{(\ell,i,j)} = \Vb^{(j,\ell)}, \qquad \gamma_\ell \ge 0, j=1,\ldots,j_*; k=1,\ldots,k^*; \ell=1,\ldots,L . \end{equation}

Consider, for example, that in general, $\mrg((\MA + \gamma I) \MU) \not= \mrg(\MA \MU)$ unless $\MU$ is an invariant subspace of $\MA$; so even something as seemly innocuous as an identity shift (i.e., $\vek E=\vek I$) can cause difficulties with any solver that falls into the recycling framework. However, these problems do have structure that can be exploited which allows for {\it extensions in
the spirit of recycling}.

In the following, we split the discussion of recycling for systems with additional structure
into scalar identity shifts, scalar non-identity shifts, and  more general continuous parameter dependence.
For scalar identity shifts, the discussion is further split into shift-dependent right-hand sides, right-hand sides independent of shifts, and multiple right-hand sides in addition to shifts.  Note that these categories are not mutually exclusive, since some methods will transform a family shift-independent
right-hand sides into one for which there is shift-dependence.  Such cases are noted in \Cref{shift-dep-rhs} with forward references to relevant
sections.

\subsection{Scalar shifted matrices $\curl{\vek A + \gamma_{\ell}\vek I}_{\ell=1}^{L}$ whose RHS changes with each shift}\label{shift-dep-rhs}
In this subsection, we consider $\ME = \MI$ for a single right-hand side and fixed system matrix in the sequence, which allows
us to drop indices $i,j$ in \eqref{eq:MRHS_all}.
We let $\vek x_{0}^{(\ell)}$ be the initial approximation for the shifted system
with shift $\gamma_\ell$, $\vek t_{j}^{(\ell)}$ to be the correction generated for this shift at iteration $j$, and
$\vek r_{j}^{(\ell)} = \vek b^{(\ell)} - (\vek A + \gamma_\ell \vek I)(\vek x_{0}^{(\ell)} + \vek t_{j}^{(\ell)})$ to be the
residual for the system associated to shift $\ell$.
A property of Krylov subspaces which makes them attractive for treating a family
of shifted systems is that for a given seed vector, the subspace is invariant with respect to scalar shifts by
the identity of the coefficient matrix.  More generally, 
\begin{equation}\label{eqn.shift-invar}
	\CK_{j}(\vek A + \gamma_{\ell_{1}}\vek I,\vek u) = \CK_{j}(\vek A + \gamma_{\ell_{2}}\vek I,\widetilde{\vek u})
\end{equation}
where $\widetilde{\vek u} = \omega\vek u$ for some $\omega\in\bbC\setminus\left\lbrace 0 \right\rbrace$ and
for any nonzero values of $\gamma_{\ell_{1}}$ and $\gamma_{\ell_{2}}$.
If $\vek b^{( \ell )} = \vek b$ 
 and
$\vek x_{0}^{(\ell)}=\vek 0$ for all $\ell$,
such as with systems arising in lattice quantum-chromodynamics \cite{Frommer1995} and Tikhonov regularization \cite{FM.1999},  
we can design (non-recycling based) solvers to take advantage of this shift invariance.
In this subsection,
we are concerned with reusing information when the right-hand sides do vary with $\gamma_\ell$.
Such shift-dependent RHS situations arise naturally in acoustics problems, but also arise when the system \eqref{eq:MRHS_all} for fixed $i,j$ denotes a
correction equation: that is, when the right-hand side denotes the initial residual $\V{r}^{(\ell)}_j $.

There are iterative methods that have been tailored to accommodate multiple shifts, and we discuss those briefly here first.  Then we move on to discussing impediments to adding recycling {\it on top of} these shifted system solvers in \Cref{sec:difficult}, and methods
that are used to overcome these obstacles in \Cref{sec.rhs-shift-indep}.

\paragraph{Iterative Methods Tailored to Shifted Systems}
For Arnoldi-based methods for large non-Hermitian
matrices, it is likely that restarting will be necessary.
Enforcing a (Petrov-) Galerkin condition for each shifted system
does not guarantee that residuals remain collinear.
At the restart stage, then, we are possibly in the position of having to solve systems with multiple shifts, and shift dependent RHS, where the RHS now correspond to the current respective residuals.
The necessary conditions for residuals to remain collinear at restart were characterized in \cite[Theorem 1]{Frommer2003}.
Methods such as the full orthogonalization method (FOM) \cite{Simoncini2003a}, conjugate gradients \cite{FM.1999,Frommer2003},
and bi-conjugate gradients (BiCG) \cite{Frommer2003} maintain a natural residual collinearity.
Other methods such as GMRES, MINRES, and QMR do not naturally maintain residual collinearity.
Since QMR and MINRES use short-term recurrences,
this was shown not to be a great obstacle \cite{Freund.Shifted-QMR.1993}.

In \cite{Frommer1998},
a restarted GMRES algorithm to simultaneously solve a family of shifted systems was derived.  The key is that the residual of
only one system is minimized.  The residuals for the other systems are then explicitly forced to be collinear to this minimized residual.
This collinearity
correction may not always exist, but it was shown in \cite{Frommer1998}
that if $\vek A$ has field of values in the right half-plane and the
shifts are all positive, real numbers, then one can always enforce the residual collinearity.

An extension of BiCGStab($\ell$) for shifted systems has also been proposed \cite{Frommer2003}.
Shifted BiCGStab($\ell$)
works by alternating $\ell$-cycles of shifted BiCG (which naturally maintains residual collinearity) and of shifted GMRES (which enforces residual collinearity).
As with shifted GMRES, shifted BiCGStab($\ell$) inherits the property that it will always be able to generate collinear residuals in the case that the shifts are all positive
and the coefficient matrix $\vek A$ has field of values in the right half-plane.

\subsubsection{Difficulties combining recycling with shifted system solvers} \label{sec:difficult}
Combining the shifted GMRES method \cite{Frommer1998} with the GCRO-DR
approach \cite{Parks.GCRODR.2006} was a part of the doctoral dissertation \cite{S.2012-thesis}, and these results were
extended and refined in \cite{SSX.2014}.  In this work, it was shown that for general non-Hermitian coefficient matrices and
augmentation subspaces, {\it it is not possible to embed a shifted restarted GMRES within a subspace recycling framework
as described in \Cref{section.framework}}.  An alternative is proposed, but its effectiveness decreases as the magnitude of the shift increases.

Essentially, for a given augmentation space $\CU$ and its image $\CC$, each shifted system must be projected as in
\Cref{section.framework}, and the shifted restarted GMRES then is applied to the projected problems.
This leads to three challenges: the projection of the initial residuals,
the exploitation of possible shift invariance of the Krylov subspace
generated by the projected coefficient matrix, and the enforcement of a residual collinearity condition at the end of
each cycle.  The correct initial residual projection was not fully treated until \cite{S.2016}, and we will defer discussion
thereof until Section \ref{sec:success.recyc.shift}.
Under favorable circumstances in the case of optimal recycling methods for minimum residual Krylov subspace methods,
it can be proven that a family of projected shifted matrices still generate the same Krylov subspace.  This fact has
been used in works as early as \cite{KilmerdeSturler2006}.
\begin{proposition}{\rm\cite[Proposition 3.1]{SSX.2014}}\label{prop.proj-shift-invar}
Let $\vek Q$ be the orthogonal projector onto $\CC$.  Then for any $\vek v\in\CC^{\perp}$ we have that
\begin{equation}\nonumber
  \CK_{j}\left(\left( \vek I - \vek Q \right)\vek A,\vek v \right)=\CK_{j}\left(\left( \vek I - \vek Q \right)\left( \vek A + \gamma\vek I \right),\vek v \right).
\end{equation}
\end{proposition}
This projected operator shift invariance was later extended to Sylvester operators in \cite{S.2016}.  What Proposition
\ref{prop.proj-shift-invar} shows is that when restarting is not required, one can generate one augmented Krylov subspace
$\CU + \CK_{j}\left(\left( \vek I - \vek Q \right)\vek A,\vek r_{0} \right)$ and compute minimum residual corrections for
each shifted system therefrom.  This was used in \cite{KilmerdeSturler2006}, since the problems were real symmetric; and as was shown,
it is compatible in the real case with complex shifts.

The question that remains is whether such a method can be extended to the non-Hermitian case when restarts are required.  Ideally,
as in \cite{Frommer1998}, for one system we would compute the minimum residual correction, and for all other systems corrections
would be computed that give residuals collinear to the minimized residual.  However, it was shown in
\cite[Theorem 1]{SSX.2014} that this is not generally possible.  However, in certain situations, the collinear residual
was shown to exist.  For example, if $\CU=\CC$ is an invariant subspace, then it is possible to compute collinear residuals as in
\cite{Frommer1998}.  More generally, if
\begin{equation}\label{eqn.collin-containment}
	\CU + \CK_{j}\left(\left( \vek I - \vek Q \right)\vek A,\vek r_{0} \right) \subset \CC + \CK_{j+1}\left(\left( \vek I - \vek Q \right)\vek A,\vek r_{0} \right),
\end{equation}
then it is possible to enforce residual collinearity.  Indeed
the shifted
GMRES-DR method \cite{Darnell2008} takes advantage of this relationship, as Morgan previously proved that a Krylov subspace
augmented with harmonic Ritz vectors satisfies \eqref{eqn.collin-containment} \cite{Morgan.GMRESDR.2002}.  Absent this property, however, enforcing shifted residuals to be collinear to the minimized one is not possible \cite{SSX.2014}.

\subsubsection{Effective recycling strategies for shifted systems}\label{sec:success.recyc.shift}
The problems discussed in \eqref{sec:difficult} are impediments to combining these two technologies, but a number of strategies
have been proposed which either overcome or work around the problems discussed in \Cref{sec:difficult}; and for scalar shifted linear systems, there are certain circumstances
for which we can take advantage of the shift invariance property and still augment.
We note again that these methods treat the case where the right-hand
sides/residuals do differ at some point. Even if we begin in the setting that we have the same right-hand side for all shifts, we are dealing with methods
which destroy that structure, leaving residuals which are either collinear but not equal or which have no relationship to one another.

\paragraph{GMRES-DR and FOM-DR} As has been mentioned, for solving a single linear system with non-Hermitian coefficient matrix using an Arnoldi-based restarting approach,
shifted versions of GMRES-DR and FOM-DR have been proposed and have shown to be effective \cite{Darnell2008}.  Both unshifted algorithms work by retaining some harmonic Ritz or Ritz
vectors, respectively, to augment the Krylov subspace generated in the next cycle.  It is shown in \cite{Morgan.GMRESDR.2002} that in each case, that the resulting space is in actuality a Krylov subspace with
a different starting vector.  In the augmentation language used in this survey, these appended vectors span the subspace $\CU$, and the space $\CU$ and $\CC = \vek A\CU$ satisfy \Cref{eqn.collin-containment}.
Thus, in the case of GMRES-DR, one can enforce residual collinearity, as in the manner of \cite{Frommer1998}.  For FOM-DR, one has the natural collinearity through enforcement of the Galerkin condition
due to \cite[Theorem 1]{Frommer2003}.

\paragraph{Direct projection} In \cite{SSX.2014}, it was shown that one cannot generally enforce shifted residuals to be collinear with the minimized residual of the seed system when minimizing over an augmented
Krylov subspace.
One can instead simply exploit the shifted system structure directly.  In \cite{S.2014.2}, it is proposed to perform a minimum residual projection for all shifted systems
over the (augmented) Krylov subspace generated by the base matrix and right-hand side.
For non-augmented Krylov subspaces, one could still exploit the shifted system structure such that
the methods was still reasonably efficient, but the version of this method for GCRO-DR is still quite costly in terms of extra floating-point calculations one must carry out for the shifted systems.

\paragraph{Sylvester equation interpretation}
It has been observed that a family of shifted systems can instead be interpreted as the Sylvester equations
\begin{equation}\nonumber
	\vek A\underbrace{\begin{bmatrix} \vek x_{1} & \vek x_{2} & \cdots \vek x_{s}	\end{bmatrix}}_{=:\vek X} + \begin{bmatrix} \vek x_{1} & \vek x_{2} & \cdots \vek x_{s} \end{bmatrix}\underbrace{{\rm diag}\curl{ \gamma_{1},  \gamma_{2}, \cdots,  \gamma_{s} } }_{=:\vek D} = \underbrace{\begin{bmatrix} \vek b^{(\gamma_{1})} & \vek b^{(\gamma_{2})} & \cdots \vek b^{(\gamma_{s})}  \end{bmatrix}}_{=:\vek B}.
\end{equation}
In \cite{S.1996}, Simoncini showed that one could approximate the solution to this Sylvester equation by generating the block Krylov subspace $\K_{j}\prn{\vek A,\vek R_{0}}$ generated by the block initial residual $\vek R_{0} = \vek B - \vek A\vek X_{0} - \vek X_{0}\vek D$.  One can apply a
GMRES minimization over this space simultaneously for all shifted systems without consideration of residual collinearity or any relationship between residuals. Thus,
building a subspace recycling method on top of this
process suffers from none of the restrictions (see \cite{S.2016}) that were seen in \cite{SSX.2014}.

This method also allows one to take advantage of theory presented in \cite{KT.2011}, where it was shown that any parameter-dependent family of linear systems
has a family of solutions well-approximated in a subspace of small dimension.  The smallness of the dimension has an upper bound depending on how smoothly the right-hand side and linear system depend on the parameter.  In \cite{S.2016}, an experiment is set up such that
a family of shifted systems (i.e., $\vek A+\gamma\vek I$ depends linearly on $\gamma$) have right-hand sides $\vek b^{(\gamma)}$ which have $\Cn{\infty}$ dependence on $\gamma$.  A few of the systems are solved using standard methods. These solutions
span a recycling subspace used to rapidly solve all other systems.  This structure can be exploited more generally in the context of recycling; cf. \Cref{section.cont-param-dep}.

\textbf{\emph{What if the residuals \underline{are} collinear?}}
If the RHS do not depend on the shifts, we can use the strategies employed in \Cref{sec.rhs-shift-indep}.  If
the residuals are collinear, one can use, e.g., the shifted GMRES method \cite{Frommer1998}, though this cannot take advantage of recycling.
Regardless, the Sylvester framework is a viable option.  If there is no recycled subspace,
one can begin by applying a cycle of GMRES to all shifted systems, generating a single Krylov subspace using the shift invariance, since residuals are collinear.
Applying the GMRES minimization to each shifted system renders the residuals non-collinear, and we are then in the general setting
and we can apply the Sylvester equation strategy.  If there is an initial recycled subspace, one applies the Sylvester equation strategy immediately,
as this will destroy residual collinearity anyway.

\subsection{Sequences of shifted systems with shift independent RHS}\label{sec.rhs-shift-indep}
In the following sections (until \Cref{section.cont-param-dep}), we discuss specific applications where the {\it coefficient matrices are real-valued}
with possible complex scalar shifting.  Therefore, Hermitian conjugation is replaced with transposition.

During the course of solving the optimization problem in image reconstruction from diffuse optical tomographic data, the authors of \cite{KilmerdeSturler2006} encounter sequences of shifted linear systems of multiple right-hand sides (MRHS)
(\ref{eq:MRHS_all}) where the RHS do not depend on the shift, and where the {\it non-zero shifts are pure imaginary\footnote{The approach can be readily modified for a real or complex shift.}} (i.e. $\gamma_\ell = \imath \frac{\omega}{\nu}$):
\begin{equation} \label{eq:MRHS} (\MA^{(i)} + \gamma_\ell \MI) \Vx^{(\ell,i,j)} = \Vb^{(j)}, \qquad \gamma_\ell \ge 0, j=1,\ldots,j_*; i=1,\ldots,i^*; \ell=1,\ldots,L, \end{equation}
where $\gamma_1 = 0$.  In their application, $\MA^{(i)}$ are real and symmetric, though it is possible to extend the idea to non-symmetric matrices.
Right-hand sides do not change as $i$ changes, and are not a function of shift\footnote{In order to generalize to non-symmetric $\MA^{(i)}$, the MINRES solver would need to be replaced with a GMRES solver.  If restarts are necessary, then it may be necessary to incorporate components from the previous section as then intermediate residuals would depend on
shift value.}.

\paragraph{No shifts, Multiple RHS} The contributions contained in the 2006 paper \cite{KilmerdeSturler2006} include recycling that takes advantage of similarities of the systems across all
three indices, tailored in part to the optimization process.
That is, for a given right-hand side, the recycle space consists of the shared recycling basis of approximate eigenvectors, augmented by a very small number of recent, prior solutions for that right-hand side, where the idea of which prior solutions to include depends on where the solve occurs in the optimization process.

\paragraph{Including Multiple Shifts} To describe their approach for the non-zero shifts we focus on a fixed right-hand side and drop the dependence on $i$ and $j$.  Let $\MU$ denote the recycle space for the current $i,j$.  Then the recycling recurrence gives
\[ \MA \MV_m = \M{C} \MB_m + \MV_{m+1} \underline{\MT}_m, \qquad \MB_m := \M{C}^T \MA \MV_{m}, \]
and as before, an optimal solution is sought in $\mrg([\MV_m,\MU])$.
After some manipulation, the least squares problem for the solution ultimately leads to
the projected problem
\[ \min_{\Vy,\Vz} \left\| \begin{bmatrix} \xi \V{e}_1 \\ \M{C}^T \Vb \\ 0 \end{bmatrix} - \begin{bmatrix}
  \underline{\MT}_m +  \gamma_\ell \underline{\MI}_m &  \gamma \MV_{m+1}^T \MU \\ \MB_m & \MI +  \gamma_\ell \M{C}^T \MU \\ 0 &  \gamma_\ell \MN \end{bmatrix} \begin{bmatrix} \Vy \\ \Vz \end{bmatrix} \right\|_2  \mbox{ or }
\min_{\Vy,\Vz} \left\| \begin{bmatrix} \M{C}^T \Vb \\ \xi \V{e}_1 \end{bmatrix} - \begin{bmatrix} \MI +  \gamma_\ell (\M{C}^T \MU) & \MB_m \\ 0 & \underline{\MT}_m +  \gamma_\ell \underline{\MI}_m \end{bmatrix}
\begin{bmatrix} \Vz \\ \Vy \end{bmatrix} \right\|_2
\]
where  $\MN$ arises from orthogonalizing $\MU$ against $[\MV_{m+1},\M{C}]$ and the solution is $\Vx = \MV_m \Vy + \MU \Vz$.
Numerical difficulties are observed if $\mrg(\MU)$ is very close to an invariant subspace of $\MA$.  The problem on the right is preferred if $\mrg(\MU)$ is expected to be very close to an invariant subspace of $\MA$, or when storage is at a premium, since, in this case,
$\MV_m \Vy$ and $\MB_m \Vy$ can be computed via short term recurrences when $\MA$ is symmetric, avoiding the storage of $\MV_m$.

\subsection{Inner-outer Recycling}
The concept of inner-outer recycling for systems of the form (\ref{eq:MRHS_all})
was first developed in
\cite{MishaMeghan.innerouter.2017}
in the context of producing a global basis necessary for producing a reduced order model (ROM) in the diffuse optical tomographic imaging problem.
Let the current global basis for creating a reduced order model be the same as the {\bf master recycle space}, $\mrg(\MU)$.
After solving the next set of systems in the sequence, columns may be appended to the global basis.
Once the global basis is sufficient, system solves for future $i$ are replaced by system solves with $\MU^T \MA^{(i)} \MU$, the ROM, which has a significantly smaller dimension than that of $\MA^{(i)}$.  Going forward, the optimization relies on solves with the ROM, so optimization steps are much cheaper.
Unfortunately for the number of columns in $\MU$ that are needed for a good ROM, orthogonalization would become too costly to use this as a recycle space in the typical setting.

The authors of \cite{MishaMeghan.innerouter.2017} observed that to enhance the global basis, one need not find an approximation to the system solutions directly, but rather to add information that is not already reconstructible from
$\mrg(\MU)$.
For each $j$, there is a correction equation
\begin{equation} \label{eq:resi}
 \MA^{(i)} \Vg^{(j)} = \Vr^{(j)}, \qquad \Vg^{(j)} := \Vx^{(j)} - \MU \M{C}^T \Vb^{(j)},
\end{equation}
(assuming fixed $i$ and $\gamma = 0$) to which the authors look for the optimal $\Vg^{(j)}$ restricted to a suitable subspace $\mathcal{S}^{(j)}$.
Their idea is to define a {\bf local recycling space} $\mrg(\MU^{(j)}) \subset \mrg(\MU)$ and use it to apply recycling MINRES to solve (\ref{eq:resi}). 
Since $\mrg(\MU^{(j)}) \subset \mrg(\MU)$ is maintained, a short-term recurrence update is also possible.


\paragraph{\bf Updating Across Shifted Systems}
The global basis matrix must provide a suitable reduced transfer function for  the 0 {\bf as well as} the pure imaginary shifts
$\gamma_\ell$.  The authors of \cite{MishaMeghan.innerouter.2017} observed that in their application, the magnitude of the shifts are small enough that the real parts of
the solutions to the shifted systems are not far (in a relative sense) from the corresponding solutions to the 0 shift case.
Thus, they initialize the master recycle space $\MU$ (and all the local recycle spaces
$\MU^{(j)}$) with the solutions to the 0 frequency systems for $k=1$.  They also augment $\MU$ (and $\MU^{(j)}$) by the {\it imaginary parts of the solutions to the corresponding imaginary shifted systems} for $k=1$.  This ensures that $\MU$ and all the $\MU^{(j)}$ remain real for all $k$.

For a $k \ge 2$, the non-shifted systems for $j = 1,\ldots,j_*$ are first solved using the inner-outer recycling described above.  Then, the non-zero shifted systems are updated.
Assuming $\Vx^{(\ell,j)} \approx \MU \Vq^{(\ell,j)}$, a 
Petrov-Galerkin projection is applied so that  $\M{C}^{T} \Vr^{(\ell,j)} = \bf 0$, giving the solution estimate for the $\ell$th shift and $j$th right-hand side as
\begin{equation}
\Vx^{(\ell,j)} \approx \MU \left( \MI + \imath \gamma_\ell \M{C}^{T} \MU \right)^{-1} \M{C}^{T} \Vb^{(j)}, \qquad \mbox{where} \qquad \M{C}^{T} \MU = \MU^{T} \MA^{(i)} \MU.
\label{eq:initXsh}
\end{equation}
No additional correction system is solved for any shifted system; the corrections are found only for the non-shifted systems for all right-hand sides, as before.
It should be noted that the elements of the work in \cite{S.2016} were inspired by this strategy.

\subsection{A family of matrices $\curl{\vek A + \gamma_{\ell}\vek E}_{\ell=1}^{L}$, $\vek E\neq\vek I$}

\paragraph{Inner-Outer recycling over shifts}
In \cite{Oconnell-thesis}, a similar idea of maintaining an outer, master subspace across all shifts, and local, shift-specific recycle spaces is introduced.
The method assumes symmetric $\MA^{(i)}$, but the shift matrix $\ME$ need only be real,
and the shifts can be anything.  We describe the procedure here for a single RHS, but details for MRHS can be found in \cite{Oconnell-thesis}.
The master subspace, $\MU$, is initially seeded with all the solutions across all the shifts for the first system ($i=1$), as well as approximate invariant subspace information for a fixed shift value (for simplicity, assume this is for $\gamma = 0$).
At step $i+1$, the columns of $\M{C}$ provide an orthonormal basis for $\mrg(\MA^{(i+1)}\MU)$.
Then, initial guesses to the solutions across all shifts are obtained over $\MU$ via a Petrov-Galerkin constraint on residuals.  If the corresponding residual $\Vr^{(\ell)}$ to the $\ell$th shift is not small enough,
$\MU^{(\ell)}$ is selected such that $\mrg(\MU^{(\ell)}) \subseteq \mrg(\MU)$.  The orthonormal columns of
$\M{C}^{(\ell)}$ provide a basis for
$\mrg( (\MA^{(i+1)} + \gamma_\ell \ME) \MU^{(\ell)})$, and recycling using this shift-specific $\MU^{(\ell)}$ and
$\MC^{(\ell)}$ is used to solve the residual correction equation.

\paragraph{Double Shifts, Single RHS}  In \cite{Saibabaetal2015}, the authors consider an application in hyperspectral DOT, but the system matrix $\MA$ remains unchanged and recycling is only over the wavelengths (i.e. shifts are real-valued).  Here,
the system matrices are a small perturbation to systems shifted by the identity,
$(\MA + \gamma_\ell \MI + \mu_\ell \ME)$ and with a single RHS.
It is also assumed $\| \ME \|_2 \ll \| \MA \|_2$.  Thus, the shift invariance property is first used to generate initial guesses to solutions, ignoring $\ME$.   Recycling is then used on each individual residual correction equation.  While the general
recycling strategy follows \cite{Parks.GCRODR.2006}, it differs in the way that the recycle and range spaces are updated.  This is possible
because $\MA_\ell := \MA + \sigma_\ell \MI + \mu_\ell \ME = \MA_{\ell -1} + \Delta_\ell$, where
$\Delta_\ell = (\sigma_\ell - \sigma_{\ell -1})\MI + (\mu_\ell - \mu_{\ell-1}) \ME$.  So, once the $\MU$ and $\M{C}$ for
$\MA$ are known, finding $\M{C}_\ell$ for $\MA_\ell$ can be done independently across $\ell$ and parallelized, so all the
shifted systems can be handled independently.

\paragraph{Other Methods}  There have been other iterative approaches proposed in the literature that deal with
shifted systems.  In some of these, such as \cite{BKLSS.mpgmres.2016}, the focus is on preconditioning to convert systems to those for which shift-invariance can be leveraged, but subspace recycling is not explicitly used and therefore we do not review that
literature here.  Other literature comes closer in spirit:  the authors of \cite{Meerbergen.bai.lanczos-param-sym.2010} propose a deflation based Lanczos
approach for the symmetric parameterized systems whose matrices correspond to a frequency response function.
The method is similar in the sense of projecting out a set of (in this case,  generalized) eigenvectors, so some form of deflation that leads to faster convergence. This can be done for multiple shifts. However, since this is based on generalized eigenvectors, this does not give standard spectral deflation in the sense described above, so we do not consider it further here.

\subsection{General parameterized families with continuous parameter dependence}\label{section.cont-param-dep}
We discussed in \Cref{sec:success.recyc.shift}, if the right-hand sides exhibit sufficiently smooth
dependence on the shift for all shifts in some interval, then the solutions associated to all shifts in this interval can be approximated
to machine accuracy in a small subspace $\CS$ which has a dimension dependent on the smoothness of the dependence on the shift.  This
theory was developed in \cite{KT.2011}, and it applies more generally to parameter-dependent linear systems of the form
\begin{equation}\label{eqn.param-dep}
\vek A(s)\vek x(s) = \vek b(s)\qquad s\in\bbC.
\end{equation}
If dependence of $\vek A(s)$ and $\vek b(s)$ is sufficiently smooth for all $s$ in some neighborhood $\Omega\subset\bbC$,
then the associated solutions all exist in
a small subspace $\CS$ of dimension $d$, which depends on the smoothness of the dependence.  In \cite{KT.2011}, upper bounds for $d$ are given
in terms of the smoothness of this dependence on $s$.  Thus, one can propose a recycling algorithm to take advantage of this theory.
For \eqref{eqn.param-dep}, if we have determined $d<d_{M}$, then we can build a recycled subspace by choosing
$\curl{s_{1},s_{2},\ldots, s_{d_{M}}}\subset\Omega$
in the neighborhood and use any method to solve the linear systems associated to these shifts.  The solutions
$\curl{\vek x(s_{1}),\vek x(s_{2}),\ldots,\vek x(s_{d_{M}})}$ form the recycled subspace which can be used to solve all systems for other $s\in\Omega$.

\section{Uses in practical applications}\label{section.applications}
\paragraph{Large Scale Software Libraries}
Recycling solvers have been implemented in major software/solver libraries, most notably in PETSc \cite{petsc-web-page,petsc-user-ref,petsc-efficient,PETScRecy,PETScDoc}, see \cite{jolivet2016block} for a discussion on performance and applications in elasticity and electromagnetics, Trilinos \cite{1089021,TrilinosTut07,BelosPckg}, and the DLR-TAU library from the German Aerospace Center, see \cite{XuTimme_17,xu2016enabling}, which also detail several challenging applications in CFD.

\paragraph{Computational Scientific and Engineering Applications}
Recycling solvers have been used in a wide range of applications, ranging from calculations for fundamental problems in computational physics to large-scale astrophysical simulations, tomography and medical imaging, and many applications in computational engineering, sometimes with modifications that serve a specific application. Recycling solvers and closely related approaches have been used in Lattice Quantum Chromodynamics \cite{aoki2010physical,aoki20092,bolten2013preconditioning,MR2599785,Darnell2008}, in particular, \cite{aoki20092} mentions that recycling is important for handling physical regimes with very small eigenvalues. Many problems in design involve sequences of slowly changing linear systems, for which recycling is highly efficient, e.g., in topology optimization and other structural optimization problems \cite{Wang.TopOptRecyc.2007,StuLeWangPaulino_06,ZhangStuShap_20,Carlberg2015,Choi-etal-WCSMO_19}, and aerodynamic shape optimization \cite{chen2019gcro,HZ.simp-FGCROT.2010}.
Recycling has also been used to compute reduced order models (for a range of
applications) \cite{Ahuja.Recyc-BiCGStab-MOD.2015,Ahuja.Recyc-BiCG.2012,feng2009parametric,feng2013subspace,MishaMeghan.innerouter.2017}.
Another important area is nonlinear optimization, such as nonlinear least-squares, e.g., in tomography \cite{MishaMeghan.innerouter.2017,KilmerdeSturler2006,Recyc-Imp-Tom.2010,Saibabaetal2015}
and blind deconvolution \cite{hennig2015probabilistic}.
Many applications arise in engineering, such as computational fluid dynamics and nonlinear structural problems
\cite{Gosselet.etal.reuse-Kryl-nonlin.2012,xu2016enabling,XuTimme_17,
Ahuja.Sturler.rGCROT-rBiCGStab-Hybrid.2015,leon2011unified,MMM.optim-DCG-elliptic-parallel.2013,Risler2000}, acoustics \cite{Meerbergen.bai.lanczos-param-sym.2010,keuchel2016combination},
and problems from electromagnetics and electrical circuits  \cite{Ye:2008:GKR:1391469.1391646,el2010field,jolivet2016block,
Peng_EMWave_11,GGPV.FGMRES-DR.2012,GirGratMart_07}.
Recycling has found many applications in uncertainty quantification and partial differential equations with stochastic components \cite{EierErnUll_07,jin2007parallel}.

\section{Outlook and future work}\label{sec.outlook}
There are yet many interesting extensions of the work mentioned above.
One important area is a better understanding of what type of subspaces to recycle
for fast convergence and how to compute such subspaces
efficiently, especially in the context of particular applications
and in terms of what can be learned from previous iterations/linear systems.
A second area is convergence theory related to various recycling approaches, particularly, biorthogonality-based recycling approaches.
Third,  further work is needed to investigate how to best combine recycling and preconditioning and to determine whether or not the the framework outlined here can help in this respect. Finally, we note that there are classes of problems -- for example, discrete ill-posed problems -- where the convergence needs and problem properties are different;  thus different ways of thinking in this context might be needed.


\section*{Acknowledgments}
The authors wish to thank the ANLA activity group for inviting us to write on this topic. The work by Eric
de Sturler was supported in part by the grant NSF DMS
1720305; the work by Misha Kilmer was supported in part by the grant
NSF DMS-1720291.  The authors would also like to thank Daniel B. Szyld for his comments on the presentation of this survey
as well as the two referees for their helpful comments which improved the quality and flow of the paper greatly.


\bibliography{Bibliography.bib}%

\begin{thebibliography}{157}
\expandafter\ifx\csname natexlab\endcsname\relax\def\natexlab#1{#1}\fi
\expandafter\ifx\csname url\endcsname\relax
  \def\url#1{{\tt #1}}\fi
\expandafter\ifx\csname urlprefix\endcsname\relax\def\urlprefix{URL }\fi
\expandafter\ifx\csname doiprefix\endcsname\relax\def\doiprefix{doi:}\fi

\bibitem[{Abdel-Rehim et~al.(2010)Abdel-Rehim, Morgan, Nicely, and
  Wilcox}]{Morgan.rest-Lanczos-defl.2010}
Abdel-Rehim, A.~M., R.~B. Morgan, D.~A. Nicely, and W.~Wilcox, 2010: Deflated
  and restarted symmetric {L}anczos methods for eigenvalues and linear
  equations with multiple right-hand sides. {\it SIAM J. Sci. Comput.\/}, {\bf
  32}, no. 1, 129--149, doi:10.1137/080727361.

\bibitem[{Abdel-Rehim et~al.(2014)Abdel-Rehim, Stathopoulos, and
  Orginos}]{EigCG-nonsymm-lanczos}
Abdel-Rehim, A.~M., A.~Stathopoulos, and K.~Orginos, 2014: Extending the
  eig{CG} algorithm to nonsymmetric {L}anczos for linear systems with multiple
  right-hand sides. {\it Numer. Linear Algebra Appl.\/}, {\bf 21}, no. 4,
  473--493, doi:10.1002/nla.1893.

\bibitem[{Agullo et~al.(2014)Agullo, Giraud, and
  Jing}]{Giraud.BGMRES-DR-Inexact.2014}
Agullo, E., L.~Giraud, and Y.-F. Jing, 2014: Block {GMRES} method with inexact
  breakdowns and deflated restarting. {\it SIAM Journal on Matrix Analysis and
  Applications\/}, {\bf 35}, no. 4, 1625--1651, doi:10.1137/140961912.

\bibitem[{Ahuja(2009)}]{A.2009}
Ahuja, K., 2009: {\it Recycling Bi-Lanczos Algorithms: Bi{CG}, {CGS}, and
  Bi{CGSTAB}\/}. Master's thesis, Virginia Polytechnical Institute, Blacksburg,
  Department of Mathematics, Blacksburg, Virginia.

\bibitem[{Ahuja et~al.(2015)Ahuja, Benner, de~Sturler, and
  Feng}]{Ahuja.Recyc-BiCGStab-MOD.2015}
Ahuja, K., P.~Benner, E.~de~Sturler, and L.~Feng, 2015: Recycling {B}i{CGSTAB}
  with an application to parametric model order reduction. {\it SIAM J. Sci.
  Comput.\/}, {\bf 37}, no. 5, S429--S446, doi:10.1137/140972433.

\bibitem[{Ahuja et~al.(2012)Ahuja, de~Sturler, Gugercin, and
  Chang}]{Ahuja.Recyc-BiCG.2012}
Ahuja, K., E.~de~Sturler, S.~Gugercin, and E.~R. Chang, 2012: Recycling
  {B}i{CG} with an application to model reduction. {\it SIAM J. Sci.
  Comput.\/}, {\bf 34}, no. 4, A1925--A1949, doi:10.1137/100801500.

\bibitem[{Al~Daas et~al.(2018)Al~Daas, Grigori, H{\'e}non, and
  Ricoux}]{al2018recycling}
Al~Daas, H., L.~Grigori, P.~H{\'e}non, and P.~Ricoux, 2018: Recycling {K}rylov
  subspaces and reducing deflation subspaces for solving sequence of linear
  systems.  INRIA Paris.

\bibitem[{Amritkar et~al.(2015)Amritkar, de~Sturler, {\'S}wirydowicz, Tafti,
  and Ahuja}]{Ahuja.Sturler.rGCROT-rBiCGStab-Hybrid.2015}
Amritkar, A., E.~de~Sturler, K.~{\'S}wirydowicz, D.~Tafti, and K.~Ahuja, 2015:
  Recycling {K}rylov subspaces for {CFD} applications and a new hybrid
  recycling solver. {\it J. Comput. Phys.\/}, {\bf 303}, 222--237,
  doi:10.1016/j.jcp.2015.09.040.

\bibitem[{Aoki et~al.(2009)Aoki, Ishikawa, Ishizuka, Izubuchi, Kadoh, Kanaya,
  Kuramashi, Namekawa, Okawa, Taniguchi, Ukawa, Ukita, Yamazaki, and
  Yoshi\'e}]{aoki20092}
Aoki, S., K.-i. Ishikawa, N.~Ishizuka, T.~Izubuchi, D.~Kadoh, K.~Kanaya,
  Y.~Kuramashi, Y.~Namekawa, M.~Okawa, Y.~Taniguchi, A.~Ukawa, N.~Ukita,
  T.~Yamazaki, and T.~Yoshi\'e, 2009: 2+1 flavor lattice {QCD} toward the
  physical point. {\it Physical review D\/}, {\bf 79}, no. 3, 034503.

\bibitem[{Aoki et~al.(2010)Aoki, Ishikawa, Ishizuka, Izubuchi, Kadoh, Kanaya,
  Kuramashi, Namekawa, Okawa, Taniguchi, Ukawa, Ukita, Yamazaki, and
  Yoshi\'e}]{aoki2010physical}
--- 2010: Physical point simulation in 2+1 flavor lattice {QCD}. {\it Physical
  Review D\/}, {\bf 81}, no. 7, 074503.

\bibitem[{Axelsson et~al.(1996)Axelsson, Neytcheva, and
  Polman}]{axelsson1996application}
Axelsson, O., M.~Neytcheva, and B.~Polman, 1996: An application of the
  bordering method to solve nearly singular systems. {\it Vestnik Moskovskogo
  Universiteta, Seria 15, Vychisl. Math. Cybern\/}, {\bf 1}, 3--25.

\bibitem[{Baglama et~al.(1998)Baglama, Calvetti, Golub, and
  Reichel}]{Bagl-etal_1998}
Baglama, J., D.~Calvetti, G.~H. Golub, and L.~Reichel, 1998: Adaptively
  preconditioned {GMRES} algorithms. {\it SIAM Journal on Scientific
  Computing\/}, {\bf 20}, no. 1, 243--269.

\bibitem[{Baglama et~al.(2003)Baglama, Calvetti, and Reichel}]{MR1978154}
Baglama, J., D.~Calvetti, and L.~Reichel, 2003: I{RBL}: an implicitly restarted
  block-{L}anczos method for large-scale {H}ermitian eigenproblems. {\it SIAM
  J. Sci. Comput.\/}, {\bf 24}, no. 5, 1650--1677,
  doi:10.1137/S1064827501397949.

\bibitem[{Baglama and Reichel(2005)}]{MR2201173}
Baglama, J. and L.~Reichel, 2005: Augmented implicitly restarted {L}anczos
  bidiagonalization methods. {\it SIAM J. Sci. Comput.\/}, {\bf 27}, no. 1,
  19--42, doi:10.1137/04060593X.

\bibitem[{Baglama and Reichel(2007{\natexlab{a}})}]{BR-2.2007}
--- 2007{\natexlab{a}}: Augmented {GMRES}-type methods. {\it Numerical Linear
  Algebra with Applications\/}, {\bf 14}, no. 4, 337--350, doi:10.1002/nla.518.

\bibitem[{Baglama and Reichel(2007{\natexlab{b}})}]{BR.2007}
--- 2007{\natexlab{b}}: Decomposition methods for large linear discrete
  ill-posed problems. {\it Journal of Computational and Applied Mathematics\/},
  {\bf 198}, no. 2, 332--343, doi:10.1016/j.cam.2005.09.025.

\bibitem[{Baglama and Reichel(2013)}]{MR3123847}
--- 2013: An implicitly restarted block {L}anczos bidiagonalization method
  using {L}eja shifts. {\it BIT\/}, {\bf 53}, no. 2, 285--310.

\bibitem[{Baglama et~al.(2013)Baglama, Reichel, and
  Richmond}]{BRR-augLSQR.2013}
Baglama, J., L.~Reichel, and D.~Richmond, 2013: An augmented {LSQR} method.
  {\it Numerical Algorithms\/}, {\bf 64}, no. 2, 263--293,
  doi:10.1007/s11075-012-9665-8.

\bibitem[{Bakhos et~al.(2017)Bakhos, Kitanidis, Ladenheim, Saibaba, and
  Szyld}]{BKLSS.mpgmres.2016}
Bakhos, T., P.~K. Kitanidis, S.~Ladenheim, A.~K. Saibaba, and D.~B. Szyld,
  2017: Multipreconditioned {GMRES} for shifted systems. {\it SIAM Journal on
  Scientific Computing\/}, {\bf 39}, no. 5, S222--S247, doi:10.1137/16M1068694.

\bibitem[{Balay et~al.(2019{\natexlab{a}})Balay, Abhyankar, Adams, Brown,
  Brune, Buschelman, Dalcin, Dener, Eijkhout, Gropp, Karpeyev, Kaushik,
  Knepley, May, McInnes, Mills, Munson, Rupp, Sanan, Smith, Zampini, Zhang, and
  Zhang}]{petsc-web-page}
Balay, S., S.~Abhyankar, M.~F. Adams, J.~Brown, P.~Brune, K.~Buschelman,
  L.~Dalcin, A.~Dener, V.~Eijkhout, W.~D. Gropp, D.~Karpeyev, D.~Kaushik, M.~G.
  Knepley, D.~A. May, L.~C. McInnes, R.~T. Mills, T.~Munson, K.~Rupp, P.~Sanan,
  B.~F. Smith, S.~Zampini, H.~Zhang, and H.~Zhang, 2019{\natexlab{a}}: {\it
  {PETS}c {W}eb page\/}. \url{https://www.mcs.anl.gov/petsc}.

\bibitem[{Balay et~al.(2019{\natexlab{b}})Balay, Abhyankar, Adams, Brown,
  Brune, Buschelman, Dalcin, Dener, Eijkhout, Gropp, Karpeyev, Kaushik,
  Knepley, May, McInnes, Mills, Munson, Rupp, Sanan, Smith, Zampini, Zhang, and
  Zhang}]{petsc-user-ref}
--- 2019{\natexlab{b}}: {PETS}c users manual.  ANL-95/11 - Revision 3.12.
  Argonne National Laboratory.
\newline\urlprefix\url{https://www.mcs.anl.gov/petsc}

\bibitem[{Balay et~al.(1997)Balay, Gropp, McInnes, and Smith}]{petsc-efficient}
Balay, S., W.~D. Gropp, L.~C. McInnes, and B.~F. Smith, 1997: Efficient
  management of parallelism in object oriented numerical software libraries.
  {\it Modern Software Tools in Scientific Computing\/}, E.~Arge, A.~M.
  Bruaset, and H.~P. Langtangen, Eds., Birkh{\"{a}}user Press, 163--202.

\bibitem[{Belos Package - Trilinos((no year))}]{BelosPckg}
Belos Package - Trilinos, (no year): {\it Belos: An iterative linear solvers
  package\/}.
\newline\urlprefix\url{https://docs.trilinos.org/dev/packages/belos/doc/html/index.html}

\bibitem[{Bolten et~al.(2013)Bolten, Bo{\v z}ovi{\'c}, and
  Frommer}]{bolten2013preconditioning}
Bolten, M., N.~Bo{\v z}ovi{\'c}, and A.~Frommer, 2013: Preconditioning of
  {K}rylov subspace methods using recycling in lattice {QCD} computations. {\it
  PAMM\/}, {\bf 13}, no. 1, 413--414.

\bibitem[{Burrage and Erhel(1998)}]{BurrErhe1998}
Burrage, K. and J.~Erhel, 1998: On the performance of various adaptive
  preconditioned {GMRES} strategies. {\it Numerical Linear Algebra with
  Applications\/}, {\bf 5}, 101--121.

\bibitem[{Calvetti et~al.(2003)Calvetti, Reichel, and Shuibi}]{CRS.2003}
Calvetti, D., L.~Reichel, and A.~Shuibi, 2003: Enriched {K}rylov subspace
  methods for ill-posed problems. {\it Linear Algebra Appl.\/}, {\bf 362},
  257--273, doi:10.1016/S0024-3795(02)00533-5.

\bibitem[{Calvetti et~al.(1994)Calvetti, Reichel, and Sorensen}]{MR1270124}
Calvetti, D., L.~Reichel, and D.~C. Sorensen, 1994: An implicitly restarted
  {L}anczos method for large symmetric eigenvalue problems. {\it Electron.
  Trans. Numer. Anal.\/}, {\bf 2}, no. March, 1--21.

\bibitem[{Carlberg et~al.(2016)Carlberg, Forstall, and Tuminaro}]{Carlberg2015}
Carlberg, K., V.~Forstall, and R.~Tuminaro, 2016: Krylov-subspace recycling via
  the {POD}-augmented conjugate-gradient method. {\it SIAM J. Matrix Anal.
  Appl.\/}, {\bf 37}, no. 3, 1304--1336.
\newline\urlprefix\url{https://doi.org/10.1137/16M1057693}

\bibitem[{Carvalho et~al.(2011)Carvalho, Gratton, Lago, and
  Vasseur}]{CGLV.FGCRODR.2012}
Carvalho, L.~M., S.~Gratton, R.~Lago, and X.~Vasseur, 2011: A flexible
  generalized conjugate residual method with inner orthogonalization and
  deflated restarting. {\it SIAM J. Matrix Anal. Appl.\/}, {\bf 32}, no. 4,
  1212--1235, doi:10.1137/100786253.

\bibitem[{Chan and Ng(1999)}]{CN.GalProjAnalMRHS.1999}
Chan, T.~F. and M.~K. Ng, 1999: Galerkin projection methods for solving
  multiple linear systems. {\it SIAM J. Sci. Comput.\/}, {\bf 21}, no. 3,
  836--850, doi:10.1137/S1064827598310227.

\bibitem[{Chan and Wan(1997)}]{Chan1997}
Chan, T.~F. and W.~L. Wan, 1997: {Analysis of projection methods for solving
  linear systems with multiple right-hand sides}. {\it SIAM Journal on
  Scientific Computing\/}, {\bf 18}, no. 6, 1698--1721.

\bibitem[{Chapman and Saad(1997)}]{Chapman1997}
Chapman, A. and Y.~Saad, 1997: Deflated and augmented {K}rylov subspace
  techniques. {\it Numerical Linear Algebra with Applications\/}, {\bf 4}, no.
  1, 43--66,
  doi:10.1002/(SICI)1099-1506(199701/02)4:1<43::AID-NLA99>3.3.CO;2-Q.

\bibitem[{Chen and Nadarajah(2019)}]{chen2019gcro}
Chen, C.-H. and S.~Nadarajah, 2019: {GCRO} with dynamic deflated restarting for
  solving adjoint systems of equations for aerodynamic shape optimization. {\it
  International Journal of Numerical Methods for Heat \& Fluid Flow\/}.

\bibitem[{Choi et~al.(2019)Choi, Oxberry, White, and
  Kirchdoerfer}]{Choi-etal-WCSMO_19}
Choi, Y., G.~Oxberry, D.~White, and T.~T. Kirchdoerfer, 2019: Accelerating
  topology optimization using reduced order models.  Lawrence Livermore
  National Lab.(LLNL), Livermore, CA (United States).

\bibitem[{Chung et~al.(2020)Chung, de~Sturler, and Jiang}]{ChuStuJia_20}
Chung, J., E.~de~Sturler, and J.~Jiang, 2020: Hybrid projection methods with
  recycling for inverse problems.  arXiv:2007.00207. arXiv Preprints.

\bibitem[{Darnell et~al.(2008)Darnell, Morgan, and Wilcox}]{Darnell2008}
Darnell, D., R.~B. Morgan, and W.~Wilcox, 2008: Deflated {GMRES} for systems
  with multiple shifts and multiple right-hand sides. {\it Linear Algebra and
  its Applications\/}, {\bf 429}, no. 10, 2415--2434,
  doi:10.1016/j.laa.2008.04.019.

\bibitem[{de~Sturler(1996)}]{EdS_HH96}
de~Sturler, E., 1996: Inner-outer methods with deflation for linear systems
  with multiple right-hand sides. {\it Householder Symposium XIII, Proceedings
  of the Householder Symposium on Numerical Algebra, Pontresina, Switzerland,
  1996.\/}, 193--196.

\bibitem[{{de Sturler}(1996)}]{EDS.GCRO.1996}
{de Sturler}, E., 1996: Nested {K}rylov methods based on {GCR}. {\it Journal of
  Computational and Applied Mathematics\/}, {\bf 67}, no. 1, 15--41,
  doi:10.1016/0377-0427(94)00123-5.

\bibitem[{{de Sturler}(1999)}]{deSturler.GCROT.1999}
--- 1999: Truncation strategies for optimal {K}rylov subspace methods. {\it
  SIAM Journal on Numerical Analysis\/}, {\bf 36}, no. 3, 864--889,
  doi:10.1137/S0036142997315950.

\bibitem[{{de Sturler}(2011)}]{desturler.recycling-convergence-abstract.2012}
--- 2011: Convergence bounds for approximate invariant subspace recycling for
  sequences of linear systems. {\it Program of the Householder Symposium XVIII
  on Numerical Linear Algebra\/}, 51--52.

\bibitem[{de~Sturler et~al.(2020)de~Sturler, Kilmer, and
  Soodhalter}]{dSKS.2018}
de~Sturler, E., M.~Kilmer, and K.~M. Soodhalter, 2020: Krylov subspace
  augmentation for the solution of shifted systems: a review. {\it in
  preparation\/}.

\bibitem[{de~Sturler et~al.(2006)de~Sturler, Le, Wang, and
  Paulino}]{StuLeWangPaulino_06}
de~Sturler, E., C.~Le, S.~Wang, and G.~H. Paulino, 2006: Large scale topology
  optimization using preconditioned {K}rylov subspace recycling and continuous
  approximation of material distribution. {\it Proceedings Multiscale and
  Functionally Graded Materials 2006 (M\&FGM 2006), Oahu Island (Hawaii),
  15--18 October 2006\/}, G.~H. Paulino, M.-J. Pindera, J.~R.~H.~Dodds, F.~A.
  Rochinha, E.~Dave, and L.~Chen, Eds., AIP, number 973 in AIP Conference
  Proceedings, 279 -- 284, iSBN: 978-0-7354-0492-2.

\bibitem[{Dolean et~al.(2015)Dolean, Jolivet, and
  Nataf}]{dolean2015introduction}
Dolean, V., P.~Jolivet, and F.~Nataf, 2015: {\it An Introduction to Domain
  Decomposition Methods: Algorithms, Theory, and Parallel Implementation\/}.
  SIAM, Philadelphia.

\bibitem[{Dong et~al.(2014)Dong, Garde, and Hansen}]{DGH.2014}
Dong, Y., H.~Garde, and P.~C. Hansen, 2014: R{${}^3$}{GMRES}: including prior
  information in {GMRES}-type methods for discrete inverse problems. {\it
  Electron. Trans. Numer. Anal.\/}, {\bf 42}, 136--146.

\bibitem[{Eiermann et~al.(2000)Eiermann, Ernst, and Schneider}]{Eiermann2000}
Eiermann, M., O.~G. Ernst, and O.~Schneider, 2000: Analysis of acceleration
  strategies for restarted minimal residual methods. {\it Journal of
  Computational and Applied Mathematics\/}, {\bf 123}, no. 1-2, 261--292,
  doi:10.1016/S0377-0427(00)00398-8, numerical analysis 2000, Vol. III. Linear
  algebra.

\bibitem[{Eiermann et~al.(2007)Eiermann, Ernst, and Ullmann}]{EierErnUll_07}
Eiermann, M., O.~G. Ernst, and E.~Ullmann, 2007: Computational aspects of the
  stochastic finite element method. {\it Comput. Vis. Sci\/}, {\bf 10}, no. 1,
  3--15, doi:https://doi.org/10.1007/s00791-006-0047-4.

\bibitem[{Eirola and Nevanlinna(1989)}]{Eirola.Nevanlinna.1989}
Eirola, T. and O.~Nevanlinna, 1989: Accelerating with rank-one updates. {\it
  Linear Algebra and its Applications\/}, {\bf 121}, 511--520,
  doi:https://doi.org/10.1016/0024-3795(89)90719-2.

\bibitem[{Eisenstat et~al.(1983)Eisenstat, Elman, and Schultz}]{EisElS83}
Eisenstat, S.~C., H.~C. Elman, and M.~H. Schultz, 1983: Variational iterative
  methods for nonsymmetric systems of linear equations. {\it SIAM Journal on
  Numerical Analysis\/}, {\bf 20}, 345--357.

\bibitem[{El-Moselhy(2010)}]{el2010field}
El-Moselhy, T.~A., 2010: {\it Field solver technologies for variation-aware
  interconnect parasitic extraction\/}. Ph.D. thesis, Massachusetts Institute
  of Technology.

\bibitem[{Engl et~al.(1996)Engl, Hanke, and Neubauer}]{EHN.1996-book}
Engl, H.~W., M.~Hanke, and A.~Neubauer, 1996: {\it Regularization of inverse
  problems\/}, volume 375 of {\it Mathematics and its Applications\/}. Kluwer
  Academic Publishers Group, Dordrecht, viii+321 pp.

\bibitem[{Erhel et~al.(1996)Erhel, Burrage, and Pohl}]{ErhelBurraPohl_96}
Erhel, J., K.~Burrage, and B.~Pohl, 1996: Restarted {GMRES} preconditioned by
  deflation. {\it Journal of Computational and Applied Mathematics\/}, {\bf
  69}, 303--318.

\bibitem[{Erhel and Guyomarc'h(2000)}]{ErhGuy_2000}
Erhel, J. and F.~Guyomarc'h, 2000: An augmented conjugate gradient method for
  solving consecutive symmetric positive definite linear systems. {\it SIAM
  Journal on Matrix Analysis and Applications\/}, {\bf 21}, no. 4, 1279 --
  1299, doi:10.1137/S0895479897330194.

\bibitem[{Erlangga and Nabben(2008)}]{Erl.Nabben.Defl-balPrec.2008}
Erlangga, Y.~A. and R.~Nabben, 2008: Deflation and balancing preconditioners
  for {K}rylov subspace methods applied to nonsymmetric matrices. {\it SIAM J.
  Matrix Anal. Appl.\/}, {\bf 30}, no. 2, 684--699, doi:10.1137/060678257.

\bibitem[{Faber et~al.(1996)Faber, Joubert, Knill, and Manteuffel}]{MR1410697}
Faber, V., W.~Joubert, E.~Knill, and T.~Manteuffel, 1996: Minimal residual
  method stronger than polynomial preconditioning. {\it SIAM J. Matrix Anal.
  Appl.\/}, {\bf 17}, no. 4, 707--729, doi:10.1137/S0895479895286748.

\bibitem[{Faber and Manteuffel(1984)}]{Faber.Manteuffel.1984}
Faber, V. and T.~Manteuffel, 1984: Necessary and sufficient conditions for the
  existence of a conjugate gradient method. {\it SIAM Journal on Numerical
  Analysis\/}, {\bf 21}, no. 2, 352--362, doi:10.1137/0721026.

\bibitem[{Fadeev and Fadeeva(1963)}]{fadeev1963computational}
Fadeev, D.~K. and V.~N. Fadeeva, 1963: {\it Computational methods of linear
  algebra\/}. W.H.Freeman \& Co, New York.

\bibitem[{Feng et~al.(2009)Feng, Benner, and Korvink}]{feng2009parametric}
Feng, L., P.~Benner, and J.~G. Korvink, 2009: Parametric model order reduction
  accelerated by subspace recycling. {\it Proceedings of the 48h IEEE
  Conference on Decision and Control (CDC) held jointly with 2009 28th Chinese
  Control Conference\/}, IEEE, 4328--4333.

\bibitem[{Feng et~al.(2013)Feng, Benner, and Korvink}]{feng2013subspace}
--- 2013: Subspace recycling accelerates the parametric macro-modeling of mems.
  {\it International Journal for Numerical Methods in Engineering\/}, {\bf 94},
  no. 1, 84--110.

\bibitem[{Fischer(1998)}]{Fischer.init-guess-from-prev.1998}
Fischer, P.~F., 1998: Projection techniques for iterative solution of
  {$A\underline x=\underline b$} with successive right-hand sides. {\it Comput.
  Methods Appl. Mech. Engrg.\/}, {\bf 163}, no. 1-4, 193--204,
  doi:10.1016/S0045-7825(98)00012-7.

\bibitem[{Fletcher(1976)}]{Fletcher-BiCG.1976}
Fletcher, R., 1976: Conjugate gradient methods for indefinite systems. {\it
  Numerical analysis ({P}roc 6th {B}iennial {D}undee {C}onf., {U}niv. {D}undee,
  {D}undee, 1975)\/}, Springer, Berlin, 73--89. Lecture Notes in Mathematics,
  Vol. 506.

\bibitem[{Frank and Vuik(2001)}]{FranVuik_2001}
Frank, J. and K.~Vuik, 2001: On the construction of deflation-based
  preconditioners. {\it SIAM Journal on Scientific Computing\/}, {\bf 23}, no.
  2, 442--462, doi:10.1137/S1064827500373231.

\bibitem[{Freund(1993{\natexlab{a}})}]{Freund.Shifted-QMR.1993}
Freund, R.~W., 1993{\natexlab{a}}: Solution of shifted linear systems by
  quasi-minimal residual iterations. {\it Numerical linear algebra ({K}ent,
  {OH}, 1992)\/}, de Gruyter, Berlin, 101--121.

\bibitem[{Freund(1993{\natexlab{b}})}]{Freund.TFQMR}
--- 1993{\natexlab{b}}: A transpose-free quasi-minimal residual algorithm for
  non-{H}ermitian linear systems. {\it SIAM J. Sci. Comput.\/}, {\bf 14}, no.
  2, 470--482, doi:10.1137/0914029.

\bibitem[{Freund and Nachtigal(1991)}]{Freund.QMR}
Freund, R.~W. and N.~M. Nachtigal, 1991: Q{MR}: a quasi-minimal residual method
  for non-{H}ermitian linear systems. {\it Numerische Mathematik\/}, {\bf 60},
  no. 3, 315--339, doi:10.1007/BF01385726.

\bibitem[{Frommer(2003)}]{Frommer2003}
Frommer, A., 2003: {${\rm BiCGStab}(l)$} for families of shifted linear
  systems. {\it Computing\/}, {\bf 70}, no. 2, 87--109,
  doi:10.1007/s00607-003-1472-6.

\bibitem[{Frommer and Gl{\"a}ssner(1998)}]{Frommer1998}
Frommer, A. and U.~Gl{\"a}ssner, 1998: Restarted {GMRES} for shifted linear
  systems. {\it SIAM Journal on Scientific Computing\/}, {\bf 19}, no. 1,
  15--26, doi:10.1137/S1064827596304563.

\bibitem[{Frommer et~al.(1995)Frommer, G\"{u}sken, Lippert, N\"{o}ckel, and
  Schilling}]{Frommer1995}
Frommer, A., S.~G\"{u}sken, T.~Lippert, B.~N\"{o}ckel, and K.~Schilling, 1995:
  {Many masses on one stroke: Economic computation of quark propagators}. {\it
  International Journal of Modern Physics C\/}, {\bf 6}, 627--638.

\bibitem[{Frommer and Maa\ss(1999)}]{FM.1999}
Frommer, A. and P.~Maa\ss, 1999: Fast {CG}-based methods for
  {T}ikhonov-{P}hillips regularization. {\it SIAM J. Sci. Comput.\/}, {\bf 20},
  no. 5, 1831--1850 (electronic), doi:10.1137/S1064827596313310.

\bibitem[{Gaul(2014)}]{Gaul.2014-phd}
Gaul, A., 2014: {\it Recycling {K}rylov subspace methods for sequences of
  linear systems: Analysis and applications\/}. Ph.D. thesis, Technischen
  Universit\"at Berlin.

\bibitem[{Gaul et~al.(2013)Gaul, Gutknecht, Liesen, and Nabben}]{GGL.2013}
Gaul, A., M.~H. Gutknecht, J.~Liesen, and R.~Nabben, 2013: A framework for
  deflated and augmented {K}rylov subspace methods. {\it SIAM Journal on Matrix
  Analysis and Applications\/}, {\bf 34}, no. 2, 495--518,
  doi:10.1137/110820713.

\bibitem[{Gaul and Schl\"{o}mer(2015)}]{Gaul.S.rMINRES.2015}
Gaul, A. and N.~Schl\"{o}mer, 2015: Preconditioned recycling {K}rylov subspace
  methods for self-adjoint problems. {\it Electron. Trans. Numer. Anal.\/},
  {\bf 44}, 522--547.

\bibitem[{Giraud et~al.(2007)Giraud, Gratton, and Martin}]{GirGratMart_07}
Giraud, L., S.~Gratton, and E.~Martin, 2007: Incremental spectral
  preconditioners for sequences of linear systems. {\it Applied Numerical
  Mathematics\/}, {\bf 57}, no. 11 - 12, 1164--1180.

\bibitem[{Giraud et~al.(2010)Giraud, Gratton, Pinel, and
  Vasseur}]{GGPV.FGMRES-DR.2012}
Giraud, L., S.~Gratton, X.~Pinel, and X.~Vasseur, 2010: Flexible {GMRES} with
  deflated restarting. {\it SIAM J. Sci. Comput.\/}, {\bf 32}, no. 4,
  1858--1878, doi:10.1137/080741847.

\bibitem[{Gosselet et~al.(2013)Gosselet, Rey, and
  Pebrel}]{Gosselet.etal.reuse-Kryl-nonlin.2012}
Gosselet, P., C.~Rey, and J.~Pebrel, 2013: Total and selective reuse of
  {K}rylov subspaces for the resolution of sequences of nonlinear structural
  problems. {\it Internat. J. Numer. Methods Engrg.\/}, {\bf 94}, no. 1,
  60--83, doi:10.1002/nme.4441.

\bibitem[{Greenbaum et~al.(1996)Greenbaum, Pt\'ak, and
  Strako\v{s}}]{Greenbaum.Any-Curve-Possible-GMRES.1996}
Greenbaum, A., V.~Pt\'ak, and Z.~Strako\v{s}, 1996: Any nonincreasing
  convergence curve is possible for {GMRES}. {\it SIAM Journal on Matrix
  Analysis and Applications\/}, {\bf 17}, 465--469.

\bibitem[{Gutknecht(2014)}]{Gutknecht.biCG-aug.2014}
Gutknecht, M.~H., 2014: Deflated and augmented {K}rylov subspace methods: a
  framework for deflated {B}i{CG} and related solvers. {\it SIAM J. Matrix
  Anal. Appl.\/}, {\bf 35}, no. 4, 1444--1466, doi:10.1137/130923087.

\bibitem[{Hanke(1995)}]{H1995-book}
Hanke, M., 1995: {\it Conjugate gradient type methods for ill-posed
  problems\/}, volume 327 of {\it Pitman Research Notes in Mathematics
  Series\/}. Longman Scientific \& Technical, Harlow, iv+134 pp.

\bibitem[{Hennig et~al.(2015)Hennig, Osborne, and
  Girolami}]{hennig2015probabilistic}
Hennig, P., M.~A. Osborne, and M.~Girolami, 2015: Probabilistic numerics and
  uncertainty in computations. {\it Proceedings of the Royal Society A:
  Mathematical, Physical and Engineering Sciences\/}, {\bf 471}, no. 2179,
  20150142.

\bibitem[{Heroux et~al.(2005)Heroux, Bartlett, Howle, Hoekstra, Hu, Kolda,
  Lehoucq, Long, Pawlowski, Phipps, Salinger, Thornquist, Tuminaro,
  Willenbring, Williams, and Stanley}]{1089021}
Heroux, M.~A., R.~A. Bartlett, V.~E. Howle, R.~J. Hoekstra, J.~J. Hu, T.~G.
  Kolda, R.~B. Lehoucq, K.~R. Long, R.~P. Pawlowski, E.~T. Phipps, A.~G.
  Salinger, H.~K. Thornquist, R.~S. Tuminaro, J.~M. Willenbring, A.~Williams,
  and K.~S. Stanley, 2005: An overview of the trilinos project. {\it ACM Trans.
  Math. Softw.\/}, {\bf 31}, no. 3, 397--423,
  doi:http://doi.acm.org/10.1145/1089014.1089021.

\bibitem[{Hestenes and Stiefel(1952)}]{Hestenes.Stiefel.CG.1952}
Hestenes, M.~R. and E.~Stiefel, 1952: Methods of conjugate gradients for
  solving linear systems. {\it Journal of Research of the National Bureau of
  Standards\/}, {\bf 49}, 409--436 (1953).

\bibitem[{Hicken and Zingg(2010)}]{HZ.simp-FGCROT.2010}
Hicken, J.~E. and D.~W. Zingg, 2010: A simplified and flexible variant of
  {GCROT} for solving nonsymmetric linear systems. {\it SIAM J. Sci.
  Comput.\/}, {\bf 32}, no. 3, 1672--1694, doi:10.1137/090754674.

\bibitem[{Jackson and Robinson(1985)}]{JR.1985}
Jackson, C.~P. and P.~C. Robinson, 1985: A numerical study of various
  algorithms related to the preconditioned conjugate gradient method. {\it
  Internat. J. Numer. Methods Engrg.\/}, {\bf 21}, no. 7, 1315--1338,
  doi:10.1002/nme.1620210711.

\bibitem[{Jin et~al.(2007)Jin, Cai, and Li}]{jin2007parallel}
Jin, C., X.-C. Cai, and C.~Li, 2007: Parallel domain decomposition methods for
  stochastic elliptic equations. {\it SIAM Journal on Scientific Computing\/},
  {\bf 29}, no. 5, 2096--2114.

\bibitem[{Jolivet and Tournier(2016)}]{jolivet2016block}
Jolivet, P. and P.-H. Tournier, 2016: Block iterative methods and recycling for
  improved scalability of linear solvers. {\it SC'16: Proceedings of the
  International Conference for High Performance Computing, Networking, Storage
  and Analysis\/}, IEEE Press, 17.
\newline\urlprefix\url{https://ieeexplore.ieee.org/abstract/document/7877095}

\bibitem[{J\"{o}nsth\"{o}vel et~al.(2012)J\"{o}nsth\"{o}vel, van Gijzen,
  MacLachlan, Vuik, and Scarpas}]{Jons-etal_2012}
J\"{o}nsth\"{o}vel, T.~B., M.~B. van Gijzen, S.~MacLachlan, C.~Vuik, and
  A.~Scarpas, 2012: Comparison of the deflated preconditioned conjugate
  gradient method and algebraic multigrid for composite materials. {\it
  Computational Mechanics\/}, {\bf 50}, no. 3, 321--333.

\bibitem[{Joubert(1994)}]{MR1261463}
Joubert, W., 1994: A robust {GMRES}-based adaptive polynomial preconditioning
  algorithm for nonsymmetric linear systems. {\it SIAM Journal on Scientific
  Computing\/}, {\bf 15}, no. 2, 427--439, doi:10.1137/0915029.

\bibitem[{Kahl and Rittich(2017)}]{KR.2012}
Kahl, K. and H.~Rittich, 2017: The deflated conjugate gradient method:
  convergence, perturbation and accuracy. {\it Linear Algebra and its
  Applications\/}, {\bf 515}, 111--129, doi:10.1016/j.laa.2016.10.027.

\bibitem[{Karchenko and Yeremin(1995)}]{KarYere1995}
Karchenko, S.~A. and A.~Y. Yeremin, 1995: Eigenvalue translation based
  preconditioners for the {GMRES}(k) method. {\it Numerical Linear Algebra with
  Applications\/}, {\bf 2}, no. 1, 51--77.

\bibitem[{Keuchel et~al.(2016)Keuchel, Biermann, and von
  Estorff}]{keuchel2016combination}
Keuchel, S., J.~Biermann, and O.~von Estorff, 2016: A combination of the fast
  multipole boundary element method and {K}rylov subspace recycling solvers.
  {\it Engineering Analysis with Boundary Elements\/}, {\bf 65}, 136--146.

\bibitem[{Kilmer and de~Sturler(2006)}]{KilmerdeSturler2006}
Kilmer, M. and E.~de~Sturler, 2006: {Recycling subspace information for diffuse
  optical tomography}. {\it SIAM J. Sci. Comput.\/}, {\bf 27}, no. 6,
  2140--2166, doi:10.1137/040610271.

\bibitem[{Kilmer et~al.(2001)Kilmer, Miller, and Rappaport}]{KMR.2001}
Kilmer, M., E.~Miller, and C.~Rappaport, 2001: Q{MR}-based projection
  techniques for the solution of non-{H}ermitian systems with multiple
  right-hand sides. {\it SIAM J. Sci. Comput.\/}, {\bf 23}, no. 3, 761--780,
  doi:10.1137/S1064827599355542.

\bibitem[{Kressner and Tobler(2011)}]{KT.2011}
Kressner, D. and C.~Tobler, 2011: Low-rank tensor {K}rylov subspace methods for
  parametrized linear systems. {\it SIAM J. Matrix Anal. Appl.\/}, {\bf 32},
  no. 4, 1288--1316, doi:10.1137/100799010.

\bibitem[{Kuroiwa and Nodera(2008)}]{KN.aug-GMRES-ip}
Kuroiwa, N. and T.~Nodera, 2008: The adaptive augmented {GMRES} method for
  solving ill-posed problems. {\it ANZIAM J.\/}, {\bf 50}, no. (C), C654--C667.

\bibitem[{Lanczos(1952)}]{Lanczos.nonsym.1952ROHTUA}
Lanczos, C., 1952: Solution of systems of linear equations by
  minimized-iterations. {\it Journal of Research of the National Bureau of
  Standards\/}, {\bf 49}, 33--53.

\bibitem[{Leon et~al.(2011)Leon, Paulino, Pereira, Menezes, and
  Lages}]{leon2011unified}
Leon, S.~E., G.~H. Paulino, A.~Pereira, I.~F. Menezes, and E.~N. Lages, 2011: A
  unified library of nonlinear solution schemes. {\it Applied Mechanics
  Reviews\/}, {\bf 64}, no. 4, 040803.

\bibitem[{Malandain et~al.(2013)Malandain, Maheu, and
  Moureau}]{MMM.optim-DCG-elliptic-parallel.2013}
Malandain, M., N.~Maheu, and V.~Moureau, 2013: Optimization of the deflated
  conjugate gradient algorithm for the solving of elliptic equations on
  massively parallel machines. {\it J. Comput. Phys.\/}, {\bf 238}, 32--47,
  doi:10.1016/j.jcp.2012.11.046.

\bibitem[{Mandel(1993)}]{Mandel.BalPrec.1993}
Mandel, J., 1993: Balancing domain decomposition. {\it Comm. Numer. Methods
  Engrg.\/}, {\bf 9}, no. 3, 233--241, doi:10.1002/cnm.1640090307.

\bibitem[{{Maya Neytcheva}(2019)}]{Neytcheva.Deflation.Talk}
{Maya Neytcheva}, 2019: {\it Deflation techniques - historical development and
  advances\/}.
  \textsc{url:}~\url{https://www.maths.tcd.ie/~ksoodha/beyonddiscrete2019/wp-content/uploads/2019/06/Neytcheva2.pdf}.

\bibitem[{Meerbergen and Bai(2009/10)}]{Meerbergen.bai.lanczos-param-sym.2010}
Meerbergen, K. and Z.~Bai, 2009/10: The {L}anczos method for parameterized
  symmetric linear systems with multiple right-hand sides. {\it SIAM J. Matrix
  Anal. Appl.\/}, {\bf 31}, no. 4, 1642--1662, doi:10.1137/08073144X.

\bibitem[{Mello et~al.(2010)Mello, de~Sturler, Paulino, and
  Silva}]{Recyc-Imp-Tom.2010}
Mello, L. s. A.~M., E.~de~Sturler, G.~H. Paulino, and E.~l. C.~N. Silva, 2010:
  Recycling {K}rylov subspaces for efficient large-scale electrical impedance
  tomography. {\it Comput. Methods Appl. Mech. Engrg.\/}, {\bf 199}, no. 49-52,
  3101--3110, doi:10.1016/j.cma.2010.06.001.

\bibitem[{Meng et~al.(2014)Meng, Zhu, and Li}]{MZL.2014}
Meng, J., P.-Y. Zhu, and H.-B. Li, 2014: A block {GCROT}(m,k) method for linear
  systems with multiple right-hand sides. {\it Journal of Computational and
  Applied Mathematics\/}, {\bf 255}, 544--554, doi:10.1016/j.cam.2013.06.014.

\bibitem[{Morgan(1995)}]{Morgan.Restarted-GMRES-eig.1995}
Morgan, R.~B., 1995: A restarted {GMRES} method augmented with eigenvectors.
  {\it SIAM Journal on Matrix Analysis and Applications\/}, {\bf 16}, no. 4,
  1154--1171, doi:10.1137/S0895479893253975.

\bibitem[{Morgan(2000)}]{Morgan2000}
--- 2000: Implicitly restarted {GMRES} and {A}rnoldi methods for nonsymmetric
  systems of equations. {\it SIAM Journal on Matrix Analysis and
  Applications\/}, {\bf 21}, no. 4, 1112--1135, doi:10.1137/S0895479897321362.

\bibitem[{Morgan(2002)}]{Morgan.GMRESDR.2002}
--- 2002: G{MRES} with deflated restarting. {\it SIAM Journal on Scientific
  Computing\/}, {\bf 24}, no. 1, 20--37, doi:10.1137/S1064827599364659.

\bibitem[{Morgan(2005)}]{Morg_2005}
--- 2005: Restarted block-{GMRES} with deflation of eigenvalues. {\it Applied
  NumericalMathematics\/}, {\bf 54}, no. 2, 222--236.

\bibitem[{Morgan et~al.(2020)Morgan, Whyte, Wilcox, and
  Yang}]{morgan2020twogrid}
Morgan, R.~B., T.~Whyte, W.~Wilcox, and Z.~Yang, 2020: Two-grid deflated
  {K}rylov methods for linear equations.  2005.03070. arXiv Preprints.

\bibitem[{Nachtigal et~al.(1992)Nachtigal, Reichel, and
  N.Trefethen}]{NachReiTref_1992}
Nachtigal, N.~M., L.~Reichel, and L.~N.Trefethen, 1992: A hybrid {GMRES}
  algorithm for nonsymmetric linear systems,. {\it SIAM J. Matrix Anal.
  Appl.\/}, {\bf 13}, no. 3, 796 -- 825.

\bibitem[{Neuenhofen and Greif(2018)}]{NeuGreif_18}
Neuenhofen, M.~P. and C.~Greif, 2018: Mstab: Stabilized induced dimension
  reduction for {K}rylov subspace recycling. {\it SIAM Journal on Scientific
  Computing\/}, {\bf 40}, no. 2, B554--B571.

\bibitem[{Nicolaides(1987)}]{Nico_87}
Nicolaides, R.~A., 1987: Deflation of conjugate gradients with applications to
  boundary value problems. {\it SIAM Journal on Numerical Analysis\/}, {\bf
  24}, no. 2, 355--365, doi:10.1137/0724027.

\bibitem[{Notay(2000)}]{Notay.2000}
Notay, Y., 2000: Flexible conjugate gradients. {\it SIAM J. Sci. Comput.\/},
  {\bf 22}, no. 4, 1444--1460, doi:10.1137/S1064827599362314.

\bibitem[{O'Connell et~al.(2017)O'Connell, Kilmer, de~Sturler, and
  Gugercin}]{MishaMeghan.innerouter.2017}
O'Connell, M., M.~E. Kilmer, E.~de~Sturler, and S.~Gugercin, 2017: Computing
  reduced order models via inner-outer {K}rylov recycling in diffuse optical
  tomography. {\it SIAM J. Sci. Comput.\/}, {\bf 39}, no. 2, B272--B297,
  doi:10.1137/16M1062880.

\bibitem[{O'Connell(2016)}]{Oconnell-thesis}
O'Connell, M.~J., 2016: {\it Advanced {T}echniques in the {C}omputation of
  {R}educed {O}rder {M}odels and {K}rylov {R}ecycling for {D}iffuse {O}ptical
  {T}omography\/}. Ph.D. thesis, Tufts University, Medford, Massachussets,
  United States of America, 115 pp.
\newline\urlprefix\url{https://dl.tufts.edu/pdfviewer/2227n161k/8c97m229d}

\bibitem[{Paige and Saunders(1975)}]{Paige1975}
Paige, C.~C. and M.~A. Saunders, 1975: Solutions of sparse indefinite systems
  of linear equations. {\it SIAM Journal on Numerical Analysis\/}, {\bf 12},
  no. 4, 617--629.

\bibitem[{Paige and Saunders(1982)}]{PS-LSQR}
--- 1982: L{SQR}: an algorithm for sparse linear equations and sparse least
  squares. {\it ACM Trans. Math. Software\/}, {\bf 8}, no. 1, 43--71,
  doi:10.1145/355984.355989.

\bibitem[{Parks et~al.(2006)Parks, de~Sturler, Mackey, Johnson, and
  Maiti}]{Parks.GCRODR.2006}
Parks, M.~L., E.~de~Sturler, G.~Mackey, D.~D. Johnson, and S.~Maiti, 2006:
  Recycling {K}rylov subspaces for sequences of linear systems. {\it SIAM
  Journal on Scientific Computing\/}, {\bf 28}, no. 5, 1651--1674,
  doi:10.1137/040607277.

\bibitem[{Parks et~al.(2013)Parks, Sampath, and Nukala}]{PSS.RCG.2016}
Parks, M.~L., R.~Sampath, and P.~K. Nukala, 2013: Efficient simulation of
  large-scale {3D} fracture networks via {K}rylov subspace recycling.
  Unpublished. Communicated by author Parks to first and second authors.

\bibitem[{Parks et~al.(2016)Parks, Soodhalter, and Szyld}]{PSS.2016}
Parks, M.~L., K.~M. Soodhalter, and D.~B. Szyld, 2016: {A block Recycled GMRES
  method with investigations into aspects of solver performance}. {\it ArXiv
  e-prints 1604.01713\/}.

\bibitem[{Peng et~al.(2011)Peng, Wang, and Lee}]{Peng_EMWave_11}
Peng, Z., X.-C. Wang, and J.-F. Lee, 2011: Integral equation based domain
  decomposition method for solving electromagnetic wave scattering from
  non-penetrable objects. {\it IEEE Transactions on Antennas and
  Propagation\/}, {\bf 59}, no. 9, 3328 -- 3338, doi:10.1109/TAP.2011.2161542.

\bibitem[{PETSc - KSPHPDDM((no year))}]{PETScRecy}
PETSc - KSPHPDDM, (no year): {\it {KSPHPDDM} solver interface with the {HPDDM}
  library\/}.
\newline\urlprefix\url{https://www.mcs.anl.gov/petsc/petsc-current/docs/manualpages/KSP/KSPHPDDM.html}

\bibitem[{PETSc Documentation((no year))}]{PETScDoc}
PETSc Documentation, (no year): {\it PETSc Manual and Examples\/}. Argonne
  National Lab.
\newline\urlprefix\url{https://www.mcs.anl.gov/petsc/petsc-current/docs/}

\bibitem[{Ramlau and Soodhalter(2020)}]{RRKMS.regrecyc.2020}
Ramlau, R. and K.~M. Soodhalter, 2020: Regularized recycling. {\it In
  Preparation\/}.

\bibitem[{Rey and Risler(1998)}]{Rey1998}
Rey, C. and F.~Risler, 1998: A {R}ayleigh-{R}itz preconditioner for the
  iterative solution to large scale nonlinear problems. {\it Numerical
  Algorithms\/}, {\bf 17}, no. 3-4, 279 -- 311, doi:10.1023/A:1016680306741.

\bibitem[{Risler and Rey(1998)}]{Risler1998}
Risler, F. and C.~Rey, 1998: On the reuse of {R}itz vectors for the solution to
  nonlinear elasticity problems by domain decomposition methods. {\it Domain
  decomposition methods, 10 ({B}oulder, {CO}, 1997)\/}, Amer. Math. Soc.,
  Providence, RI, volume 218 of {\it Contemp. Math.\/}, 334--340.

\bibitem[{Risler and Rey(2000)}]{Risler2000}
--- 2000: Iterative accelerating algorithms with {K}rylov subspaces for the
  solution to large-scale nonlinear problems. {\it Numerical Algorithms\/},
  {\bf 23}, no. 1, 1--30, doi:10.1023/A:1019187614377.

\bibitem[{Saad(1987)}]{Saad-seed}
Saad, Y., 1987: On the {L}\'{a}nczos method for solving symmetric linear
  systems with several right-hand sides. {\it Math. Comp.\/}, {\bf 48}, no.
  178, 651--662, doi:10.2307/2007834.

\bibitem[{Saad(1997)}]{Saad.Deflated-Aug-Krylov.1997}
--- 1997: Analysis of augmented {K}rylov subspace methods. {\it SIAM Journal on
  Matrix Analysis and Applications\/}, {\bf 18}, no. 2, 435--449,
  doi:10.1137/S0895479895294289.

\bibitem[{Saad(2003)}]{Saad.Iter.Meth.Sparse.2003}
--- 2003: {\it Iterative Methods for Sparse Linear Systems\/}. 2ndnd ed., SIAM,
  Philadelphia.

\bibitem[{Saad and Schultz(1986)}]{Saad.GMRES.1986}
Saad, Y. and M.~H. Schultz, 1986: {GMRES}: A generalized minimal residual
  algorithm for solving nonsymmetric linear systems. {\it SIAM Journal on
  Scientific and Statistical Computing\/}, {\bf 7}, 856--869.

\bibitem[{Saad et~al.(2000)Saad, Yeung, Erhel, and Guyomarc'h}]{SYEG.2000}
Saad, Y., M.~Yeung, J.~Erhel, and F.~Guyomarc'h, 2000: A deflated version of
  the conjugate gradient algorithm. {\it SIAM Journal on Scientific
  Computing\/}, {\bf 21}, no. 5, 1909--1926, doi:10.1137/S1064829598339761.

\bibitem[{Saibaba et~al.(2015)Saibaba, Kilmer, Miller, and
  Fantini}]{Saibabaetal2015}
Saibaba, A., M.~Kilmer, E.~Miller, and S.~Fantini, 2015: {Fast algorithms for
  hyperspectral diffuse optical tomography}. {\it SIAM J. Sci. Comput.\/}, {\bf
  37}, no. 5, B712--B743.

\bibitem[{Sala et~al.(2007)Sala, Heroux, and Day}]{TrilinosTut07}
Sala, M., M.~A. Heroux, and D.~M. Day, 2007: {\it Trilinos Tutorial\/}. Sandial
  National Laboratories.
\newline\urlprefix\url{https://trilinos.github.io/pdfs/Trilinos8.0Tutorial.pdf}

\bibitem[{Sheikh et~al.(2016)Sheikh, Lahaye, Ramos, Nabben, and
  Vuik}]{Shei-etal_2016}
Sheikh, A.~H., D.~Lahaye, L.~G. Ramos, R.~Nabben, and C.~Vuik, 2016:
  Accelerating the shifted {L}aplace preconditioner for the {H}elmholtz
  equation by multilevel deflation. {\it Journal of Computational Physics\/},
  {\bf 322}, 473--490.

\bibitem[{Sifuentes et~al.(2013)Sifuentes, Embree, and Morgan}]{Sifuentes2011}
Sifuentes, J.~A., M.~Embree, and R.~B. Morgan, 2013: G{MRES} convergence for
  perturbed coefficient matrices, with application to approximate deflation
  preconditioning. {\it SIAM J. Matrix Anal. Appl.\/}, {\bf 34}, no. 3,
  1066--1088, doi:10.1137/120884328.

\bibitem[{Simoncini(1996)}]{S.1996}
Simoncini, V., 1996: On the numerical solution of {$AX-XB=C$}. {\it BIT\/},
  {\bf 36}, no. 4, 814--830, doi:10.1007/BF01733793.

\bibitem[{Simoncini(2000)}]{Simoncini2000}
--- 2000: On the convergence of restarted {Krylov} subspace methods. {\it SIAM
  Journal on Matrix Analysis and Applications\/}, {\bf 22}, no. 2, 430--452.

\bibitem[{Simoncini(2003)}]{Simoncini2003a}
--- 2003: Restarted full orthogonalization method for shifted linear systems.
  {\it BIT. Numerical Mathematics\/}, {\bf 43}, no. 2, 459--466,
  doi:10.1023/A:1026000105893.

\bibitem[{Simoncini and Gallopoulos(1995)}]{SimonGallo_1995}
Simoncini, V. and E.~Gallopoulos, 1995: An iterative method for nonsymmetric
  systems with multiple right hand sides. {\it SIAM Journal on Scientific
  Computing\/}, {\bf 16}, no. 4, 917--933.

\bibitem[{Simoncini and Szyld(2005)}]{Simoncini2005}
Simoncini, V. and D.~B. Szyld, 2005: On the occurrence of superlinear
  convergence of exact and inexact {K}rylov subspace methods. {\it SIAM
  Review\/}, {\bf 47}, no. 2, 247--272, doi:10.1137/S0036144503424439.

\bibitem[{Simoncini and Szyld(2010)}]{Simoncini2010a}
--- 2010: Interpreting {IDR} as a {P}etrov-{G}alerkin method. {\it SIAM Journal
  on Scientific Computing\/}, {\bf 32}, 1898--1912, doi:10.1137/070685804.

\bibitem[{Sleijpen and Fokkema(1993)}]{bicgstabl.1993}
Sleijpen, G. L.~G. and D.~R. Fokkema, 1993: Bi{CG}stab{$(l)$} for linear
  equations involving unsymmetric matrices with complex spectrum. {\it
  Electron. Trans. Numer. Anal.\/}, {\bf 1}, no. Sept., 11--32 (electronic
  only).

\bibitem[{Sonneveld(1989)}]{Sonneveld.CGS.1989}
Sonneveld, P., 1989: C{GS}, a fast {L}anczos-type solver for nonsymmetric
  linear systems. {\it SIAM J. Sci. Statist. Comput.\/}, {\bf 10}, no. 1,
  36--52, doi:10.1137/0910004.

\bibitem[{Sonneveld and van Gijzen(2008)}]{Sonneveld2008}
Sonneveld, P. and M.~B. van Gijzen, 2008: { IDR($s$)}: a family of simple and
  fast algorithms for solving large nonsymmetric systems of linear equations.
  {\it SIAM Journal on Scientific Computing\/}, {\bf 31}, no. 2, 1035--1062,
  doi:10.1137/070685804.

\bibitem[{Soodhalter(2012)}]{S.2012-thesis}
Soodhalter, K.~M., 2012: {\it Krylov subspace methods with fixed memory
  requirements: Nearly Hermitian linear systems and subspace recycling\/}.
  Ph.D. thesis, Temple University.

\bibitem[{Soodhalter(2016{\natexlab{a}})}]{S.2016}
--- 2016{\natexlab{a}}: Block {Krylov} subspace recycling for shifted systems
  with unrelated right-hand sides. {\it SIAM Journal on Scientific
  Computing\/}, {\bf 38}, no. 1, A302--A324, doi:10.1137/140998214.

\bibitem[{Soodhalter(2016{\natexlab{b}})}]{S.2014.2}
--- 2016{\natexlab{b}}: {Two recursive GMRES-type methods for shifted linear
  systems with general preconditioning}. {\it Electronic Transactions on
  Numerical Analysis\/}, {\bf 45}, 499--523.
\newline\urlprefix\url{http://arxiv.org/abs/1403.4428}

\bibitem[{Soodhalter(2020)}]{KMS-aug-arn-tikh.2020}
--- 2020: Augmented {A}rnoldi-{T}ikhonov methods for ill-posed problems. {\it
  In Preparation\/}.

\bibitem[{Soodhalter et~al.(2014)Soodhalter, Szyld, and Xue}]{SSX.2014}
Soodhalter, K.~M., D.~B. Szyld, and F.~Xue, 2014: Krylov subspace recycling for
  sequences of shifted linear systems. {\it Applied Numerical Mathematics\/},
  {\bf 81}, 105--118.
\newline\urlprefix\url{http://www.sciencedirect.com/science/article/pii/S0168927414000208}

\bibitem[{Stathopoulos and Orginos(2010)}]{MR2599785}
Stathopoulos, A. and K.~Orginos, 2010: Computing and deflating eigenvalues
  while solving multiple right-hand side linear systems with an application to
  quantum chromodynamics. {\it SIAM J. Sci. Comput.\/}, {\bf 32}, no. 1,
  439--462, doi:10.1137/080725532.

\bibitem[{van~der Sluis and van~der Vorst(1986)}]{SluiVors_86}
van~der Sluis, A. and H.~van~der Vorst, 1986: The rate of convergence of
  {C}onjugate {G}radients. {\it Numerische Mathematik\/}, {\bf 48}, 543--560.

\bibitem[{van~der Vorst(1992)}]{vandervorst.bicgstab.1992}
van~der Vorst, H.~A., 1992: Bi-{CGSTAB}: a fast and smoothly converging variant
  of {B}i-{CG} for the solution of nonsymmetric linear systems. {\it SIAM
  Journal on Scientific and Statistical Computing\/}, {\bf 13}, no. 2,
  631--644, doi:10.1137/0913035.

\bibitem[{van~der Vorst(2009)}]{vanDerVorst-book.2009}
--- 2009: {\it Iterative {K}rylov methods for large linear systems\/},
  volume~13 of {\it Cambridge Monographs on Applied and Computational
  Mathematics\/}. Cambridge University Press, Cambridge, xiv+221 pp., reprint
  of the 2003 original.

\bibitem[{{van der Vorst} and Vuik(1994)}]{VV.GMRESr.1994}
{van der Vorst}, H.~A. and K.~Vuik, 1994: G{MRESR}: a family of nested {GMRES}
  methods. {\it Numerical Linear Algebra with Applications\/}, {\bf 1}, no. 4,
  369--386, doi:10.1002/nla.1680010404.

\bibitem[{Wang et~al.(2007)Wang, de~Sturler, and
  Paulino}]{Wang.TopOptRecyc.2007}
Wang, S., E.~de~Sturler, and G.~H. Paulino, 2007: Large-scale topology
  optimization using preconditioned {K}rylov subspace methods with recycling.
  {\it International Journal for Numerical Methods in Engineering\/}, {\bf 69},
  no. 12, 2441--2468, doi:10.1002/nme.1798.

\bibitem[{Xu and Timme(2017)}]{XuTimme_17}
Xu, S. and S.~Timme, 2017: Robust and efficient adjoint solver for complex flow
  conditions? {\it Computers and Fluids\/}, {\bf 148}, 26--38,
  doi:http://dx.doi.org/10.1016/j.compfluid.2017.02.012.

\bibitem[{Xu et~al.(2016)Xu, Timme, and Badcock}]{xu2016enabling}
Xu, S., S.~Timme, and K.~J. Badcock, 2016: Enabling off-design linearised
  aerodynamics analysis using {K}rylov subspace recycling technique. {\it
  Computers \& Fluids\/}, {\bf 140}, 385--396.

\bibitem[{Ye et~al.(2008)Ye, Zhu, and Phillips}]{Ye:2008:GKR:1391469.1391646}
Ye, Z., Z.~Zhu, and J.~R. Phillips, 2008: Generalized {Krylov} recycling
  methods for solution of multiple related linear equation systems in
  electromagnetic analysis. {\it Proceedings of the 45th Annual Design
  Automation Conference\/}, ACM, New York, NY, USA, DAC '08, 682--687.

\bibitem[{Zhang et~al.(2020)Zhang, de~Sturler, and Shapiro}]{ZhangStuShap_20}
Zhang, X.~S., E.~de~Sturler, and A.~Shapiro, 2020: Topology optimization with
  many right hand sides using a mirror descent stochastic approximation -
  reduction from many to a single sample. {\it Journal of Applied Mechanics\/},
  {\bf 87}, 051005--1, doi:https://doi.org/10.1115/1.4045902.

\end{thebibliography}

%
%
%

\end{document}